\documentclass[11pt]{article}
\usepackage{graphicx}
\usepackage{amsmath}
\usepackage{amsthm}
\usepackage{amssymb}
\usepackage[mathscr]{eucal}

\title{Problems on automorphism groups of nonpositively curved polyhedral complexes and
their lattices\\
{\it \normalsize To Bob Zimmer on his $60^{\rm th}$ birthday}}
\author{Benson Farb, Chris Hruska and Anne Thomas \thanks{The first author is supported in part by the NSF.
The second author is supported by the NSF under Grant Nos.\ DMS-0505659
and DMS-0731759.  The third author is supported by the NSF under Grant No.\ DMS-0805206.}}

\newtheorem{theorem}{Theorem}

\newtheorem{prob}[theorem]{Problem}
\newtheorem{conjecture}[theorem]{Conjecture}

\newtheorem{ques}[theorem]{Question}

\def\example{{\bf {\bigskip}{\noindent}Example: }}
\def\examples{{\bf {\bigskip}{\noindent}Examples: }}

\def\title{\em}

\DeclareMathOperator\Lk{Lk} \DeclareMathOperator\rank{rank}

\DeclareMathOperator\Id{Id} \DeclareMathOperator\SL{SL}
\DeclareMathOperator\PSL{PSL}
\DeclareMathOperator\SO{SO} \DeclareMathOperator\GL{GL}
\DeclareMathOperator\Sp{Sp} \DeclareMathOperator\SU{SU}
 
\DeclareMathOperator\Aut{Aut} 
 \DeclareMathOperator\Comm{Comm}
 \DeclareMathOperator\Ad{Ad}
\DeclareMathOperator\St{St} 

\DeclareMathOperator\ch{char}

\def\polhk#1{\setbox0=\hbox{#1}{\ooalign{\hidewidth
    \lower1.0ex\hbox{$\,\lhook$}\hidewidth\crcr\unhbox0}}}
\newcommand{\Swiatkowski}{\'Swi{\polhk{a}}tkowski}

\newcommand\R{\mbox{\bf R}}

\newcommand\hyp{\mbox{\bf H}}
\newcommand\Z{\mbox{\bf Z}}
\newcommand\Q{\mbox{\bf Q}}
\newcommand\C{\mbox{\bf C}}
\renewcommand\to{\longrightarrow}
\newcommand\F{\mbox{\bf F}}
\newcommand\X{\mbox{\bf X}}
\newcommand\G{\Gamma}
\newcommand{\abs}[1]{\left| {#1} \right|}

\begin{document}

\maketitle \tableofcontents

\section{Introduction}

The goal of this paper is to present a number of problems about automorphism groups of nonpositively curved
polyhedral complexes and their lattices. This topic lies at the juncture of two slightly different cultures. In
geometric group theory, universal covers of $2$--complexes are studied as geometric and topological models of
their fundamental groups, and an important way of understanding groups is to construct ``nice" actions on cell
complexes, such as cubical complexes.  From a different perspective, automorphism groups of connected,
simply-connected, locally-finite simplicial complexes may be viewed as locally compact topological groups, to
which we can hope to extend the theory of algebraic groups and their discrete subgroups. In the classification of
locally compact topological groups, these automorphism groups are natural next examples to study after algebraic
groups.  In this paper we pose some problems meant to highlight possible directions for future research.

Let $G$ be a locally-compact topological group with left-invariant Haar measure $\mu$.  A {\em lattice} (resp.
uniform lattice) in $G$ is a discrete subgroup $\Gamma<G$ with $\mu(\Gamma\backslash G)<\infty$ (resp.
$\G\backslash G$ compact). The classical study of Lie groups and their lattices was extended to algebraic groups
$G$ over nonarchimedean local fields $K$ by Ihara, Bruhat--Tits, Serre and many others. This was done by realizing
$G$ as a group of automorphisms of the Bruhat--Tits (Euclidean) building $X_G$, which is a
$\rank_K(G)$--dimensional, nonpositively curved (in an appropriate sense) simplicial complex. More recently,
Kac--Moody groups $G$ have been studied by considering the action of $G$ on the associated (twin) Tits buildings
(see, for example, Carbone--Garland~\cite{CG} and R\'emy--Ronan~\cite{RR}).

The simplest example in the algebraic case is $G=\SL(n,K)$, where one can take $K=\Q_p$ (where $\ch(K)=0$) or
$K=\F_p((t))$ (where $\ch(K)=p>0$). When $n=2$, i.e. $\rank_K(G)=1$, the building $X_G$ is the regular simplicial
tree of degree $p+1$. One can then extend this point of view to study the full group of simplicial automorphisms
of a locally-finite tree as a locally-compact topological group, and investigate the properties of the lattices it
contains.  This leads to the remarkably rich theory of ``tree lattices'', to which we refer the reader to the book
of Bass--Lubotzky \cite{BL1} as the standard reference.

One would like to build an analogous theory in dimensions $2$ and higher, with groups like $\SL(n,\Q_p)$ and
$\SL\bigl(n,\F_p((t))\bigr)$, for $n \geq 3$, being the ``classical examples''.   The increase in dimension makes life much
harder, and greatly increases the variety of phenomena that occur.

Now, let $X$ be a locally-finite, connected, simply-connected simplicial complex.  The group $G=\Aut(X)$ of
simplicial automorphisms of $X$ naturally has the structure of a locally-compact topological group, where a
decreasing neighborhood basis of the identity consists of automorphisms of $X$ which are the identity on bigger
and bigger balls. With the right normalization of the Haar measure $\mu$, due to Serre~\cite{S}, there is a useful
combinatorial formula for the covolume of a discrete subgroup $\G <G$: $$\mu(\Gamma\backslash G)=\sum_{v\in
A}\frac{1}{|\Gamma_v|}$$ where the sum is taken over vertices $v$ in a fundamental domain $A\subseteq X$ for the
$\Gamma$--action, and $|\Gamma_v|$ is the order of the $\Gamma$--stabilizer of $v$.  A discrete subgroup $\G$ is a
{\em lattice} if and only if this sum converges, and $\G$ is a {\em uniform lattice} if and only if the fundamental domain $A$
is compact.

In this paper we concentrate on the case when $\dim(X)=2$.  Most questions also make
sense in higher dimensions, where even less is understood.  When $X$ is a product of trees much is known (see, for
example, Burger--Mozes~\cite{BM}). However, the availability of projections to trees makes this a special (but
deep) theory; we henceforth assume also that $X$ is not a product. There are several themes we wish to explore,
many informed by the classical (algebraic) case and the theory of tree lattices in \cite{BL1}.   We also hope
that classical cases may be re-understood from a new, more geometric point of view.  Part of our inspiration for
this paper came from Lubotzky's beautiful paper \cite{Lu}, where he discusses the theory of tree lattices in
relation to the classical (real and $p$--adic) cases.

This paper is not meant to be encyclopedic.  It is presenting a list of problems from a specific
and biased point of view.  An important criterion in our choice of problem is that it presents
some new phenomenon, or requires some new technique or viewpoint in order to solve it. After some
background in Section~\ref{s:background}, we describe the main known
examples of polyhedral complexes
and their lattices in Section~\ref{s:examples}.  We have grouped problems on the structure of the
complex $X$ itself together with basic group-theoretic and topological properties of $\Aut(X)$ in
Section~\ref{s:properties_X}.  Section~\ref{s:comparisons} focusses on whether important
properties of linear groups and their lattices hold in this setting, while Section~\ref{s:lattices}
discusses group-theoretic properties of lattices in $\Aut(X)$ themselves.

We would like to thank Noel Brady and John Crisp for permission to use Figure~\ref{fig:BradyCrisp}, and Laurent Saloff-Coste for helpful
discussions.  We would also like to thank Fr\'{e}d\'{e}ric Haglund for making many
useful comments which greatly improved the exposition of this paper.

\section{Some background}\label{s:background}

This preliminary material is mostly drawn from Bridson--Haefliger~\cite{BH}.  We give the key definitions for
polyhedral complexes in Section~\ref{ss:polyhedral}.  (Examples of polyhedral complexes are described in
Section~\ref{s:examples} below.)  Conditions for a polyhedral complex $X$ to have nonpositive curvature, and some
the consequences for $X$, are recalled in Section~\ref{ss:npc}.  The theory of complexes of groups, which is used to construct both polyhedral complexes and their lattices, is sketched in Section~\ref{ss:cxs_of_groups}.

\subsection{Polyhedral complexes}\label{ss:polyhedral}

Polyhedral complexes may be viewed as generalizations of (geometric realizations of) simplicial complexes.  The
quotient of a simplicial complex by a group acting by simplicial automorphisms is not necessarily simplicial,
and so we work in this larger category.  Roughly speaking, a polyhedral complex is obtained by gluing together
polyhedra from some constant curvature space by isometries along faces.

More formally, let $\X^n$ be $S^n$, $\R^n$ or $\hyp^n$, endowed with Riemannian metrics of constant curvature 1,
0 and $-1$ respectively. A \emph{polyhedral complex} $X$ is a finite-dimensional CW--complex such that: \begin{enumerate} \item each
open cell of dimension $n$ is isometric to the interior of a compact convex polyhedron in $\X^n$; and \item for
each cell $\sigma$ of $X$, the restriction of the attaching map to each open codimension one face of $\sigma$ is an
isometry onto an open cell of $X$. \end{enumerate} A polyhedral complex is said to be (piecewise)
\emph{spherical}, \emph{Euclidean} or \emph{hyperbolic} if $\X^n$ is $S^n$, $\R^n$ or $\hyp^n$ respectively.
Polyhedral complexes are usually not thought of as embedded in any space.  A $2$--dimensional polyhedral complex
is called a \emph{polygonal complex}.

Given a polyhedral complex $X$, we write $G = \Aut(X)$ for the group of automorphisms, or cellular isometries, of
$X$.  A subgroup $H \leq G$ is said to act \emph{without inversions} on $X$ if for every cell $\sigma$ of $X$, the
setwise stabilizer of $\sigma$ in $H$ is equal to its pointwise stabilizer.  Note that any subgroup $H \leq G$
acts without inversions on the barycentric subdivision of $X$. The quotient of a polyhedral complex by a group
acting without inversions is also a polyhedral complex so that the quotient map is a local isometry.

Let $x$ be a vertex of an $n$--dimensional polyhedral complex $X$. The \emph{link} of $x$, written $\Lk(x,X)$, is
the spherical polyhedral complex obtained by intersecting $X$ with an $n$--sphere of sufficiently small radius
centered at $x$.  For example, if $X$ has dimension $2$, then $\Lk(x,X)$ may be identified with the graph having
vertices the $1$--cells of $X$ containing $x$ and edges the $2$--cells of $X$ containing $x$; two vertices in the
link are joined by an edge in the link if the corresponding $1$--cells in $X$ are contained in a common $2$--cell. The link may also
be thought of as the space of directions, or of germs of geodesics, at the vertex $x$.  By rescaling so that for
each $x$ the $n$--sphere around $x$ has radius say $1$, we induce a metric on each link, and we may then speak of
isometry classes of links of $X$.

\subsection{Nonpositive curvature}\label{ss:npc}

In this section, we recall conditions under which the metrics on the cells of $X$, a Euclidean or hyperbolic
polyhedral complex, may be pieced together to obtain a global metric which is respectively CAT$(0)$ or CAT$(-1)$.
Some of the consequences for $X$ are then described.

Any polyhedral complex $X$ has an intrinsic pseudometric $d$, where for $x,y \in X$, the value of
$d(x,y)$ is the infimum of lengths of paths $\Sigma$ from $x$ to $y$ in $X$, such that the restriction of
$\Sigma$ to each cell of $X$ is geodesic.  Bridson~\cite{BH} showed that if $X$ has only finitely many
isometry types of cells, for example if $G=\Aut(X)$ acts cocompactly, then $(X,d)$ is a complete
geodesic metric space.

Now assume $X$ is a Euclidean (respectively, hyperbolic) polyhedral complex such that $(X,d)$ is a complete
geodesic space.  By the Cartan--Hadamard Theorem, if $X$ is locally CAT$(0)$ (respectively, locally CAT$(-1)$), then the
universal cover $\widetilde{X}$ is CAT$(0)$ (respectively, CAT$(-1)$).  Thus to see whether a simply-connected $X$
has a CAT$(0)$ metric, we need only check a neighborhood of each point $x \in X$.

If $\dim(X) = n$ and $x$ is in the interior of an $n$--cell of $X$, then a neighborhood of $x$ is isometric to a
neighborhood in Euclidean (respectively, hyperbolic) $n$--space.  If $x$ is not a vertex but is in the
intersection of two $n$--cells, then it is not hard to see that a neighborhood of $x$ is also CAT$(0)$
(respectively, CAT$(-1)$).  Hence, the condition that $X$ be CAT$(0)$ comes down to a condition on
the nieghborhoods of the vertices of $X$, that is on their links.

There are two special cases in which it is easy to check whether neighborhoods of vertices are CAT$(0)$ or CAT$(-1)$.
These are when $\dim(X)=2$, and when $X$ is a cubical complex (defined below, and discussed in
Section~\ref{ss:CAT0_cube}).

\begin{theorem}[Gromov Link Condition]\label{t:gromov_link} A $2$--dimensional Euclidean (respectively, hyperbolic) polyhedral complex
$X$ is locally CAT$(0)$ (respectively, CAT$(-1)$) if and only if for every vertex $x$ of $X$, every injective loop in the
graph $\Lk(x,X)$ has length at least $2\pi$. \end{theorem}

Let $I^n = [0,1]^n$ be the cube in $\R^n$ with edge lengths 1.  A \emph{cubical complex} is a Euclidean polyhedral complex
with  all $n$--cells isometric to $I^n$.  Let $L$ be a simplicial complex.  We say $L$ is a \emph{flag} complex if
whenever $L$ contains the $1$--skeleton of a simplex, it contains the simplex (``no empty triangles").

\begin{theorem}[Gromov]\label{t:flag} A finite-dimensional cubical complex $X$ is locally CAT$(0)$ if and only if the link $L$ of
each vertex of $X$ is a flag simplicial complex. \end{theorem}

In general, let $X$ be a polyhedral complex of piecewise constant curvature $\kappa$ (so $\kappa = 0$ for
$X$ Euclidean, and $\kappa = -1$ for $X$ hyperbolic).

\begin{theorem}[Gromov]\label{t:CAT1} If $X$ is a polyhedral complex of piecewise constant curvature $\kappa$,
and $X$ has finitely many isometry types of cells, then $X$ is
locally \textup{CAT($\kappa$)} if and only if for all
vertices $x$ of $X$, the link $\Lk(x,X)$ is a \textup{CAT$(1)$} space. \end{theorem}

The condition that a metric space be nonpositively curved has a number of implications, described
for example in~\cite{BH}.  We highlight the
following results:
\begin{itemize} \item Any CAT$(0)$ space $X$ is contractible.
\item Let $X$ be a
complete CAT$(0)$ space.  If a group $\Gamma$ acts by isometries on $X$ with a bounded orbit, then $\Gamma$ has a
fixed point in $X$. \end{itemize} In particular, suppose $X$ is a locally finite CAT$(0)$ polyhedral complex and
$\Gamma<\Aut(X)$ is a finite group acting on $X$.  Then $\Gamma$ is contained in the stabilizer of some cell of $X$.

\subsection{Complexes of groups}\label{ss:cxs_of_groups}

The theory of complexes of groups, due to Gersten--Stallings~\cite{St} and Haefliger~\cite{Hae,BH}, generalizes
Bass--Serre theory to higher dimensions.  It may be used to construct both polyhedral complexes and lattices in
their automorphism groups.  We give here only the main ideas and some examples, and refer the reader to~\cite{BH}
for further details.

Throughout this section, if $Y$ is a polyhedral complex, then $Y'$ will denote the first barycentric subdivision of $Y$.
This is a simplicial complex with vertices $V(Y')$ and edges $E(Y')$.  Each $a \in E(Y')$ corresponds to cells
$\tau \subset \sigma$ of $Y$, and so may be oriented from $i(a)=\sigma$ to $t(a)=\tau$.  Two edges $a$ and $b$ of
$Y'$ are \emph{composable} if $i(a) = t(b)$, in which case there exists an edge $c=ab$ of $Y'$ such that $i(c) =
i(b)$, $t(c) = t(a)$ and $a$, $b$ and $c$ form the boundary of a $2$--simplex in $Y'$.

A \emph{complex of groups} $G(Y)=(G_\sigma, \psi_a, g_{a,b})$ over a polyhedral complex $Y$ is given by:
\begin{enumerate} \item a group $G_\sigma$ for each $\sigma \in V(Y')$, called the \emph{local group} at $\sigma$;
\item a monomorphism $\psi_a\colon G_{i(a)}\rightarrow G_{t(a)}$ for each $a \in E(Y')$; and \item for each pair of
composable edges $a$, $b$ in $Y'$, an element $g_{a,b} \in G_{t(a)}$, such that \[ \Ad(g_{a,b})\circ\psi_{ab} =
\psi_a \circ\psi_b \] where $\Ad(g_{a,b})$ is conjugation by $g_{a,b}$ in $G_{t(a)}$, and for each triple of
composable edges $a,b,c$ the following cocycle condition holds \[\psi_a(g_{b,c})\,g_{a,bc} = g_{a,b}\,g_{ab,c}\]
\end{enumerate}
If all $g_{a,b}$ are trivial, the complex of groups is \emph{simple}.  To date, most applications
have used only simple complexes of groups. In the case $Y$ is $2$--dimensional, the local groups of a complex of
groups over $Y$ are often referred to as face, edge and vertex groups.

\example Let $P$ be a regular right-angled hyperbolic $p$--gon, $p \geq 5$, and let $q$ be a positive integer
$\geq 2$.  Let $G(P)$ be the following \emph{polygon of groups} over $P$.  The face group is trivial, and each
edge group is the cyclic group $\Z / q\Z$.  The vertex groups are the direct products of adjacent edge groups.
All monomorphisms are natural inclusions, and all $g_{a,b}$ are trivial.

\vspace{3mm}

Let $G$ be a group acting without inversions on a polyhedral complex $X$.  The action of $G$ induces a complex of
groups, as follows. Let $Y = G \backslash X$ with $p\colon X \rightarrow Y$  the natural projection. For each $\sigma
\in V(Y')$, choose $\tilde\sigma \in V(X')$ such that $p(\tilde\sigma) = \sigma$. The local group $G_\sigma$ is
the stabilizer of $\tilde\sigma$ in $G$, and the $\psi_a$ and $g_{a,b}$ are defined using further choices. The
resulting complex of groups $G(Y)$ is unique (up to isomorphism).

Let $G(Y)$ be a complex of groups.  Then one defines the {\em fundamental group of $G(Y)$},
denoted by $\pi_1\bigl(G(Y)\bigr)$, as well as the universal cover of $G(Y)$, denoted
by $\widetilde{G(Y)}$, and an action of $\pi_1(G(Y))$ without inversion
on $\widetilde{G(Y)}$.  The quotient of $\widetilde{G(Y)}$ by this action is naturally isomorphic to $Y$, and for each cell $\sigma$ of $Y$ the stabilizer of any lift $\widetilde{\sigma}\subset\widetilde{G(Y)}$
is a homomorphic image of $G_\sigma$.  The complex of groups is called {\em developable}
whenever each homomorphism $G_\sigma\to {\rm Stab}_{\pi_1(G(Y))}(\sigma)$ is injective.
Equivalently, a complex of groups is developable if it is isomorphic to the complex of groups associated as above to an action without inversion on a simply-connected polyhedral complex.

Unlike graphs of groups, complexes of groups are not in general developable:

\example (K. Brown)  Let $G(Y)$ be the triangle of groups with trivial face group and edge groups infinite
cyclic, generated by say $a$, $b$, and $c$.  Each vertex group is isomorphic to the Baumslag--Solitar group
$BS(1,2)=\langle \,x, y \mid xyx^{-1} = y^2\, \rangle$, where the generators $x$ and $y$ are identified with the
generators of the adjacent edge groups.  The fundamental group of $G(Y)$ then has presentation $\langle\, a, b, c
\mid aba^{-1} = b^2, bcb^{-1} = c^2, cac^{-1} = a^2\,\rangle$.  It is an exercise that this is the trivial group.
Thus $G(Y)$ is not developable.

\vspace{3mm}

We now describe a local condition for developability.  Let $Y$ be a connected polyhedral complex
and let $\sigma \in
V(Y')$.  The \emph{star} of $\sigma$, written $\St(\sigma)$, is the union of the interiors of the simplices in
$Y'$ which meet $\sigma$.  If $G(Y)$ is a complex of groups over $Y$ then, even if $G(Y)$ is not developable, each
$\sigma \in V(Y')$ has a \emph{local development}. That is, we may associate to $\sigma$ an action of $G_\sigma$
on the star $\St(\tilde\sigma)$ of a vertex $\tilde\sigma$ in some simplicial complex, such that $\St(\sigma)$ is
the quotient of $\St(\tilde\sigma)$ by the action of $G_\sigma$.  To determine the local development, its link may
be computed in combinatorial fashion.

\example Suppose $G(Y)$ is a simple polygon of groups, with $G_\sigma=V$ a vertex group, with adjacent edge groups
$E_1$ and $E_2$, and with face group $F$.  We identify the groups $E_1$, $E_2$ and $F$ with their images in $V$.  The
link $L$ of the local development at $\sigma$ is then a bipartite graph.  The two sets of vertices of $L$
correspond to the cosets of $E_1$ and $E_2$ respectively in $V$, and the edges of $L$ correspond to cosets of $F$
in $V$.  The number of edges between vertices $g_1 E_1$ and $g_2 E_2$ is equal to the number of cosets of $F$ in
the intersection $g_1 E_1 \cap g_2 E_2$.  In the polygon of groups $G(P)$ given above, the link of
the local development at each vertex of $P$ will be the complete bipartite graph $K_{q,q}$.

\vspace{3mm}

If $G(Y)$ is developable, then for each $\sigma \in V(Y')$, the local development  $\St(\tilde\sigma)$ is
isomorphic to the star of each lift $\tilde\sigma$ of $\sigma$ in the universal cover $\widetilde{G(Y)}$. The
local development has a metric structure induced by that of the polyhedral complex $Y$.  We say that a complex of
groups $G(Y)$ is \emph{nonpositively curved} if for all $\sigma \in V(Y')$, the star $\St(\tilde\sigma)$ is CAT$(0)$ in this induced metric. The importance of this condition is given by:

\begin{theorem}[Stallings \cite{St}, Haefliger \cite{Hae,BH}]\label{t:nonpos} A nonpositively curved complex of groups is developable.
\end{theorem}

\example The polygon of groups $G(P)$ above is nonpositively curved and thus developable.  The links are the complete bipartite graph $K_{q,q}$
with edge lengths $\frac{\pi}{2}$, and so Gromov's Link Condition (Theorem~\ref{t:gromov_link} above) is
satisfied.

\vspace{3mm}

Let $G(Y)$ be a developable complex of groups, with universal cover a locally finite polyhedral complex $X$, and
fundamental group $\G$.  We say that $G(Y)$ is \emph{faithful} if the action of $\G$ on $X$ is faithful.  If so,
$\G$ may be regarded as a subgroup of $\Aut(X)$.  Moreover, $\Gamma$ is discrete if and only if all local groups
of $G(Y)$ are finite, and $\G$ is a uniform lattice if and only if $Y$ is compact.

\example Let $G(P)$ be the (developable) polygon of groups above, with fundamental group say $\Gamma$ and
universal cover say $X$.  Then $G(P)$ is faithful since its face group is trivial.  As all the local groups are
finite, and $P$ is compact, $\Gamma$ may be identified with a uniform lattice in $\Aut(X)$.

\section{Examples of polyhedral complexes and their lattices}\label{s:examples}

In this section we present the most studied
examples of locally finite polyhedral complexes $X$ and their lattices.  For
each case, we give the key definitions, and sketch known constructions of $X$ and of lattices in $\Aut(X)$.  There
is some overlap between examples, which we describe.  We will also try to indicate the distinctive flavor of each
class.  While results on existence of $X$ are recalled here, we defer questions of uniqueness of $X$, given
certain local data, to Section~\ref{ss:localdata} below.  Existence of lattices is also discussed further, in
Section~\ref{ss:existence_classification}.

Many of the examples we discuss are buildings, which form an important class of nonpositively curved polyhedral
complexes.  Roughly, buildings may be thought of as highly symmetric complexes, which contain many flats, and
often have algebraic structure.  Classical buildings are those associated to groups such as $SL(n,\Q_p)$, and play
a similar role for these groups to that of symmetric spaces for real Lie groups.  The basic references
for buildings are Ronan~\cite{Ron} and Brown~\cite{Br}.  A much more comprehensive treatment by
Abramenko--Brown~\cite{AB} is to appear shortly.  These works adopt a combinatorial approach.  For our purposes we
present a more topological definition, from~\cite{HP2}.

Recall that a \emph{Coxeter group} is a group $W$ with a finite generating set $S$ and presentation of the form
\[
   W = \bigl\langle\, s \in S \bigm| (s_i s_j)^{m_{ij}} = 1\,\bigr\rangle
\]
where $s_i,s_j \in S$, $m_{ii} = 1$, and if $i \not =
j$ then $m_{ij}$ is an integer $\geq 2$ or $m_{ij} = \infty$, meaning that there is no relation between $s_i$ and
$s_j$.  The pair $(W,S)$, or $(W,I)$ where $I$ is the finite indexing set of $S$, is called a \emph{Coxeter system}.  A spherical, Euclidean or hyperbolic \emph{Coxeter
polytope} of dimension $n$ is an $n$--dimensional compact convex polyhedron $P$ in the appropriate space, with
every dihedral angle of the form $\pi/m$ for some integer $m \geq 2$ (not necessarily the same $m$ for
each angle).  The group $W$ generated by reflections in the codimension one faces of a Coxeter polytope $P$ is a Coxeter group, and its action generates
a tesselation of the space by copies of $P$.

\begin{dfinition}\label{d:building} Let $P$ be an $n$--dimensional spherical, Euclidean or hyperbolic Coxeter
polytope.  Let $W = (W,S)$ be the Coxeter group generated by the set of reflections $S$ in the codimension one
faces of $P$. A \emph{spherical, Euclidean or hyperbolic building} of type $(W,S)$ is a polyhedral complex $X$
equipped with a maximal family of subcomplexes, called \emph{apartments}, each polyhedrally isometric to the
tesselation of respectively $S^n$, $\R^n$ or $\hyp^n$ by the images of $P$ under $W$ (called \emph{chambers}),
such that: \begin{enumerate} \item any two chambers of $X$ are contained in a common apartment; and \item for any
two apartments $A$ and $A'$, there exists a polyhedral isometry from $A$ onto $A'$ which fixes $A \cap A'$.
\end{enumerate}  \end{dfinition}

The links of vertices of $n$--dimensional buildings are spherical buildings of dimension $n-1$, with the induced
apartment and chamber structure.  Using this, and Theorem~\ref{t:CAT1} above, it follows that Euclidean
(respectively, hyperbolic) buildings are CAT$(0)$ (respectively,
$\textup{CAT}(-1)$).   Since buildings are
such important examples in the theory, we will spend some time describing them in detail.

\subsection{Euclidean buildings}\label{ss:Eucl_buildings}

Euclidean buildings are also sometimes known as {\em affine buildings}, or {\em
buildings of affine type}.  A simplicial tree
$X$ is a $1$--dimensional Euclidean building of type $(W,S)$, where $W$ is the infinite dihedral group, acting on
the real line with fundamental domain $P$ an interval.  The chambers of $X$ are the edges of the tree, and the
apartments $X$ are the geodesic lines in the tree.  Since the product of two buildings is also a building, it follows
that products of trees are higher-dimensional (reducible) Euclidean buildings (see Section~\ref{ss:product_trees}
below).  In this section we consider Euclidean buildings $X$ of dimension $n \geq 2$ which are not products.

\subsubsection{Classical Euclidean buildings}\label{sss:classical}

Classical Euclidean buildings are those Euclidean buildings which are associated to algebraic groups, as we now
outline.  We first construct the building for $G=\SL(n,K)$ where $K$ is a nonarchimedean local field, in terms of
lattices in $K^n$ and then in terms of $BN$--pairs (defined below).  We then indicate how the latter construction
generalizes to other algebraic groups.  Our treatment is based upon~\cite{Br}.

Let $K$ be a field.  We recall that a \emph{discrete valuation} on $K$ is a surjective homomorphism
$v\colon K^* \to \Z$, where $K^*$ is the multiplicative group of nonzero elements of $K$, such that
\[
   v(x + y) \geq \min \bigl\{v(x),v(y)\bigr\}
\]
for all $x, y \in K^*$ with $x + y \not = 0$.  We set $v(0) = +\infty$, so that $v$ is defined
and the above inequality holds for all of $K$.  A discrete valuation induces an absolute value $\abs{x} = e^{-v(x)}$
on $K$, which satisfies the nonarchimedean inequality
\[
   \abs{x + y} \leq \max\bigl\{ \abs{x}, \abs{y}\bigr\}
\]
A metric on $K$ is
obtained by setting $d(x,y) = \abs{x - y}$.
The set
$\mathcal{O} = \bigl\{ \,x \in K \bigm| \abs{x} \leq 1 \,\bigr\}$
is a subring of $K$
called the \emph{ring of integers}.  The ring $\mathcal{O}$ is compact and open in the metric topology induced by
$v$.  Pick an element $\pi \in K$ with $v(\pi) = 1$, called a \emph{uniformizer}. Every $x \in K^*$ is then
uniquely expressible in the form $x = \pi^n u$ where $n \in \Z$ and $u$ is a unit of $\mathcal{O}^*$ (so $v(u) =
0$). The ideal $\pi \mathcal{O}$ generated by $\pi$ is a maximal ideal, since every element of $\mathcal{O}$ not
in $\pi \mathcal{O}$ is a unit. Hence $k = \mathcal{O}/\pi \mathcal{O}$ is a field, called the \emph{residue field}.

\examples
\begin{enumerate}
\item For any prime $p$ the $p$--adic valuation $v$ on the field of rationals
$\Q$ is
defined by $v(x) = n$, where $x=p^n a / b$ and $a$ and $b$ are integers not divisible by $p$. The field of
$p$--adics $K = \Q_p$ is the completion of $\Q$ with respect to the metric induced by $v$, and the valuation $v$
extends to $\Q_p$ by continuity. The ring of integers is the ring of $p$--adic integers $\Z_p$, and we may take
$\pi=p$ as uniformizer. The residue field of $\Q_p$ is then the finite field
$k=\F_p$.
\item Let $q$ be a power of a prime $p$.
The field $K = \F_q((t))$ of formal Laurent series with coefficients in the finite field $\F_q$ has valuation $v$
given by
\[
   v \Bigg( \sum_{j = -m}^\infty a_j t^j \Bigg) = -m
\]
where $a_{-m}\neq 0$,
a uniformizer is $t$, and the ring of integers is
the ring of formal power series $\F_q[[t]]$.  The residue field is $k=\F_q$. \end{enumerate} \vspace{3mm}

A \emph{local nonarchimedean field} is a field $K$ which is complete with respect to the metric induced by a
discrete valuation, and whose residue field is finite.  Examples are $K=\Q_p$, which has $\ch(K)=0$, and
$K=\F_q((t))$, which has $\ch(K) = p > 0$.  In fact, all local nonarchimedean fields arise as finite extensions of
these examples.

We now fix $K$ to be a local nonarchimedean field, $\mathcal{O}$ its ring of integers, $\pi$ a uniformizer, and $k$ its
residue field. The {\em Euclidean building} associated to the group $G = \SL(n,K)$ is the geometric realization
$\abs{\Delta}$ of the abstract simplicial complex $\Delta$ which we now describe.

Let $V$ be the vector space $K^n$.  A \emph{lattice} in $V$ is an $\mathcal{O}$--submodule $L \subset V$ of the
form $L = \mathcal{O}v_1 \oplus \cdots \oplus\mathcal{O}v_n$ for some basis $\{v_1,\ldots,v_n \}$ of $V$.  If $L'$
is another lattice, then we may choose a basis $\{v_1,\ldots,v_n\}$ for $L$ such that $L'$ admits the basis $\{
\lambda_1 v_1, \ldots, \lambda_n v_n \}$ for some $\lambda_i \in K^*$.  The $\lambda_i$ may be taken to be powers
of $\pi$. Two lattices $L$ and $L'$ are \emph{equivalent} if $L = \lambda L'$ for some $\lambda \in K^*$. We write
$[L]$ for the equivalence class of $L$, and $[v_1,\ldots,v_n]$ for the equivalence class of the lattice with basis
$\{v_1,\ldots,v_n\}$.

The abstract simplicial complex $\Delta$ is defined to have vertices the set of equivalence classes of lattices in
$V$.  To describe the higher-dimensional simplices of $\Delta$, we introduce the following incidence relation. (An
\emph{incidence relation} is a relation which is reflexive and symmetric.) Two equivalence classes of lattices
$\Lambda$ and $\Lambda'$ are \emph{incident} if they have representatives $L$ and $L'$ such that \[ \pi L \subset
L' \subset L\] This relation is symmetric, since $\pi L' \in \Lambda'$ and $\pi L \in \Lambda$ satisfy \[ \pi L'
\subset \pi L \subset L'\] The simplices of $\Delta$ are then defined to be the finite sets of pairwise incident
equivalence classes of lattices in $V$.

By the definition of incidence, every top-dimensional simplex of $\Delta$ has vertex set \[[v_1, \ldots,v_i,\pi
v_{i+1}, \ldots, \pi v_n] \mbox{ for }i = 1,\ldots ,n,\] for some basis $\{v_1,\ldots,v_n\}$ of $V$. Hence
$\Delta$ is a simplicial complex of dimension $n - 1$. The geometric realization $X=\abs{\Delta}$ is thus a Euclidean polyhedral complex of dimension $n - 1$.  We note
that $n-1$ is equal to the $K$--rank of $G=\SL(n,K)$.

We now construct a simplicial complex isomorphic to $\Delta$, using certain subgroups $B$ and $N$ of $G=\SL(n,K)$.
For now, we state without proof that, with the correct Euclidean metrization, $X=|\Delta|$ is indeed a Euclidean
building, with chambers its $(n-1)$--cells, and that the vertex set of an apartment of $X$ is the set of
equivalence classes \[[\pi^{m_1}v_1, \ldots, \pi^{m_n}v_n]\] where the $m_i$ are integers $\geq 0$, and
$\{v_1,\ldots,v_n\}$ is a fixed basis for $V$.

Observe that the group $G=\SL(n,K)$ acts on the set of lattices in $V$. This action preserves equivalence of
lattices and the incidence relation, so $G$ acts without inversions on $X$.  Let $\{e_1, \ldots,e_n\}$ be the
standard basis of $V$.  We define the \emph{fundamental chamber} of $X$ to be the simplex with vertices \[
[e_1,\ldots,e_i,\pi e_{i+1},\ldots,\pi e_n], \mbox{ for }i = 1,\ldots, n,\] and the \emph{fundamental
apartment} of $X$ to be the subcomplex with vertex set \[[\pi^{m_1}e_1, \ldots, \pi^{m_n}e_n], \mbox{ where $m_i
\geq 0$}\]

Define $B$ to be the stabilizer in $G$ of the fundamental chamber, and $N$ to be  the stabilizer in $G$ of the
fundamental apartment.  There is a surjection $\SL(n,\mathcal{O}) \to \SL(n,k)$ induced by the surjection
$\mathcal{O} \to k$.  It is not hard to verify that $B$ is the inverse image in $\SL(n,\mathcal{O})$ of the upper
triangular subgroup of $\SL(n,k)$, and that $N$ is the monomial subgroup of $\SL(n,K)$ (that is, the set of
matrices with exactly one nonzero entry in each row and each column).  We say that a subgroup of $G$ is
\emph{special} if it contains a coset of $B$.

Now, from the set of cosets in $G$ of special subgroups, we form a partially ordered set, ordered by opposite
inclusion. There is an abstract simplicial complex $\Delta(G,B)$ associated to this poset. The vertices of
$\Delta(G,B)$ are cosets of special subgroups, and the simplices of $\Delta(G,B)$ correspond to chains of opposite
inclusions. Using the action of $G$ on $\Delta$ and the construction of $\Delta(G,B)$, it is not hard to see that
$\Delta(G,B)$ is isomorphic to (the barycentric subdivision of) $\Delta$.

We now generalize the construction of $\Delta(G,B)$ to algebraic groups besides $G=SL(n,K)$. Let $G$ be an
absolutely almost simple, simply connected linear algebraic group defined over $K$. Examples other than $\SL(n,K)$
include $\Sp(2n,K)$, $\SO(n,K)$, and $\SU(n,K)$.  All such groups $G$ have a Euclidean $BN$--pair, which we now define.
A \emph{$BN$--pair} is a pair of subgroups $B$ and $N$ of $G$, such that: \begin{itemize} \item $B$ and
$N$ generate $G$; \item the subgroup $T = B \cap N$ is normal in $N$; and \item the quotient $W = N/T$ admits a
set of generators $S$ satisfying certain (technical) axioms, which ensure that $(W,S)$ is a Coxeter
system.\end{itemize} A $BN$--pair is \emph{Euclidean} if the group $W$ is a Euclidean Coxeter group.  The letter $B$ stands for the
Borel subgroup, $T$ for the torus, $N$ for the normalizer of the torus, and $W$ for the Weyl group.

For $G=\SL(n,K)$, the $B$ and $N$ defined above, as $G$--stabilizers of the fundamental chamber and fundamental
apartment of $\Delta$, are a $BN$--pair. Their intersection $T$ is the diagonal subgroup of
$\SL(n,\mathcal{O})$.  The group $W$ acts on the fundamental apartment of $\Delta$ with quotient the fundamental
chamber, and is in fact isomorphic to the Coxeter group generated by reflections in the codimension one faces of a
Euclidean $(n-1)$--simplex (with certain dihedral angles).

For any group $G$ with a Euclidean $BN$--pair, one may construct the simplicial complex $\Delta(G,B)$ from
the poset of cosets of special subgroups, as described above.  The geometric realization of $\Delta(G,B)$
is a Euclidean polyhedral complex, of dimension equal to the $K$--rank of $G$.  To prove that the geometric
realization of $\Delta(G,B)$ is a building, one uses the axioms for a $BN$--pair, results about Coxeter
groups, and the Bruhat--Tits decomposition of $G$.

For classical Euclidean buildings $X$, there is a close relationship between the algebraic group $G$ to which this
building is associated, and the group $\Aut(X)$, so long as $\dim(X) \geq 2$.

\begin{theorem}[Tits~\cite{T3}]\label{t:G_AutX} Let $G$ be an absolutely almost simple, simply-connected linear algebraic group
defined over a nonarchimedean local field $K$.  Let $X$ be the Euclidean building for $G$.  If $\rank_K(G) \geq
2$, then $G$ has finite index in $\Aut(X)$ when $\ch(K) = 0$, and is cocompact in $\Aut(X)$ when $\ch(K) = p > 0$.
\end{theorem}

Thus the lattice theory of $\Aut(X)$ is very similar to that of $G$.  Existence and construction of lattices in
groups $G$ as in Theorem~\ref{t:G_AutX} are well-understood.  If $\ch(K)=0$ then $G$ does not have a nonuniform lattice
(Tamagawa~\cite{Ta}), but does admit a uniform lattice, constructed by arithmetic means
(Borel--Harder~\cite{BoHa}).  If $\ch(K) = p > 0$ then $G$ has an arithmetic nonuniform lattice, and an arithmetic
uniform lattice if and only if $G=SL(n,K)$ (Borel--Harder~\cite{BoHa}).  In real rank at least $2$ (for example, if
$G=SL(n,K)$, for $n \geq 3$) every lattice of $G$ is arithmetic (Margulis~\cite{M}).

\subsubsection{Nonclassical Euclidean buildings}\label{sss:nonclassical}

Nonclassical Euclidean buildings are those Euclidean buildings (see Definition~\ref{d:building}) which are not
the building for any algebraic group $G$ over a nonarchimedean local field.  Tits constructed uncountably many isometry classes of
nonclassical Euclidean buildings~\cite{T4}.  Nonclassical buildings may also be constructed as universal covers of
finite complexes, a method developed by Ballmann--Brin~\cite{BB1}, and examples of this kind were obtained by
Barr\'e~\cite{Ba} as well.  Ronan~\cite{Ronan86} used a construction similar to the inductive construction of Ballmann--Brin, described in
Section~\ref{ss:kL_complexes} below, to construct $2$--dimensional nonclassical Euclidean buildings.

Very few lattices are known for nonclassical buildings. In~\cite{CMSZ}, exotic lattices which act simply
transitively on the vertices of various classical and nonclassical Euclidean buildings (of type $\tilde{A}_2$) are
constructed by combinatorial methods.

\subsection{Products of trees}\label{ss:product_trees}

Let $T_1$ and $T_2$ be locally finite simplicial trees.  The product space $T_1 \times T_2$ is a polygonal
complex, where each $2$--cell is a square (edge $\times$ edge), and the link at each vertex is a complete
bipartite graph.  Products of more than two trees may also be studied.

The group $G=\Aut(T_1 \times T_2)$ is isomorphic to $\Aut(T_1) \times \Aut(T_2)$ (with a semidirect product with
$\Z/2\Z$ if $T_1 = T_2$). Thus any subgroup of $G$ may be projected to the factors. Because of this availability
of projections, the theory of lattices for products of trees is a special (but deep) theory.  See, for example,
the work of Burger--Mozes~\cite{BM}.  Many of the problems listed below may be posed in this context, but we omit
questions specific to this case.

\subsection{Hyperbolic buildings}\label{ss:hyp_buildings}

The simplest example of a hyperbolic building is \emph{Bourdon's building} $I_{p,q}$, defined and studied
in~\cite{B1}.  Here $p$ and $q$ are integers, $p \geq 5$ and $q \geq 2$.  The building $I_{p,q}$ is the (unique)
hyperbolic polygonal complex such that each $2$--cell (chamber) is isometric to a regular right-angled hyperbolic
$p$--gon $P$, and the link at each vertex is the complete bipartite graph $K_{q,q}$. The apartments of $I_{p,q}$
are hyperbolic planes tesselated by copies of $P$.  Bourdon's building is
CAT$(-1)$, and may be regarded as a
hyperbolic version of the product of two $q$--regular trees, since it has the same links.  However, $I_{p,q}$ is
not globally a product space.  The example of a polygon of groups $G(P)$ given in Section~\ref{ss:cxs_of_groups}
above has universal cover $I_{p,q}$, and the fundamental group $\Gamma$ of this polygon of groups is a uniform
lattice in $\Aut(I_{p,q})$.

Bourdon's building is a \emph{Fuchsian building}, that is, a hyperbolic building of dimension two.  More general
Fuchsian buildings have all chambers hyperbolic $k$--gons, $k \geq 3$, with each vertex angle of the form
$\pi/m$, for some integer $m \geq 2$ (depending on the vertex).  The link at each vertex with angle
$\pi/m$ is a one-dimensional spherical building $L$ which is a
generalized $m$--gon, that is, a graph with diameter $m$ and girth $2m$. For example, a complete bipartite graph is
a generalized $2$--gon.

Unlike Euclidean buildings, hyperbolic buildings do not exist in arbitrary dimension.  This is because there is a
bound ($n \leq 29$), due to Vinberg~\cite{Vin}, on the dimension $n$ of a compact convex hyperbolic Coxeter
polytope.  Gaboriau--Paulin~\cite{GP} broadened the definition of building given above (Definition~\ref{d:building}) to
allow hyperbolic buildings with noncompact chambers, in which case there are examples in any dimension, with
chambers for example ideal hyperbolic simplexes.

Various constructions of hyperbolic buildings are known.  In low dimensions, right-angled buildings (see
Section~\ref{ss:rab}) may be equipped with the structure of a hyperbolic building.  In particular, Bourdon's
building is a right-angled building.  Certain hyperbolic buildings arise as Kac--Moody
buildings (see Section~\ref{ss:KM_buildings} below), and some Davis--Moussong complexes may also be metrized as hyperbolic buildings
(see Section~\ref{ss:DM_complexes}).  Vdovina constructed some Fuchsian buildings as universal covers of finite complexes~\cite{V}. Fuchsian buildings
were constructed as universal covers of polygons of groups by Bourdon~\cite{B1,B2} and by
Gaboriau--Paulin~\cite{GP}. Haglund--Paulin~\cite{HP2} have constructed $3$--dimensional hyperbolic buildings
using ``tree-like" decompositions of the corresponding Coxeter systems.

Many of these constructions of hyperbolic buildings $X$ also yield lattices in $\Aut(X)$.  When a hyperbolic
building $X$ is a Kac--Moody building then a few lattices in $\Aut(X)$ are known from Kac--Moody theory (see for
example~\cite{R1}), and when $X$ is a Davis--Moussong complex for a Coxeter group $W$ then $W$ may be regarded as
a uniform lattice in $\Aut(X)$.  If $X$ is the universal cover of a finite complex, the fundamental group of that
complex is a uniform lattice in $\Aut(X)$. As described in Section~\ref{ss:cxs_of_groups}, if $X$ is the universal
cover of a finite complex of finite groups, such as a polygon of finite groups, then the fundamental group of the complex of groups is a uniform
lattice in $\Aut(X)$.  More elaborate complexes of groups were used by Thomas to construct both uniform and
nonuniform lattices for certain Fuchsian buildings in~\cite{Th3}.  In~\cite{B2} Bourdon was able to ``lift" lattices
for affine buildings to uniform and nonuniform lattices for certain Fuchsian buildings.

\subsection{Right-angled buildings}\label{ss:rab}

Recall that $(W,I)$ is a \emph{right-angled} Coxeter system if all the $m_{ij}$ with $i \not = j$ equal $2$ or
$\infty$.  A building $X$ of type $(W,I)$ is then a \emph{right-angled building}.  Products of trees are examples
of right-angled buildings, with associated Coxeter group the direct product of infinite dihedral groups.

Bourdon's building $I_{p,q}$, discussed in Section~\ref{ss:hyp_buildings} above, is another basic example of a
right-angled building. The Coxeter group $W$ here is generated by reflections in the sides of a regular
right-angled hyperbolic $p$--gon.  Right-angled Coxeter polytopes exist only in dimensions $n \leq 4$, and this
bound is sharp (Potyagailo--Vinberg~\cite{PV}).  Thus right-angled buildings may be metrized as hyperbolic
buildings (with compact chambers) only in dimensions $\leq 4$.

We may broaden the definition of building given above (Definition~\ref{d:building}) to allow apartments which are
Davis--Moussong complexes for $W$ (see Section~\ref{ss:DM_complexes} below), rather than just the manifold $S^n$,
$\R^n$ or $\hyp^n$ tesselated by the action of $W$.  With this definition, Gromov-hyperbolic right-angled buildings, equipped with a
piecewise Euclidean metric, exist in arbitrary dimensions (Januszkiewicz--\Swiatkowski~\cite{JS}).

The following construction of a right-angled building $X$ and a uniform lattice in $\Aut(X)$
appears in~\cite{HP1}; this construction was previously known to Davis and Meier.  It is a generalization of the polygon of groups $G(P)$ in
Section~\ref{ss:cxs_of_groups} above.  Let $(W,I)$ be a right-angled Coxeter system, and
$\{q_i\}_{i \in I}$ a set of cardinalities with $q_i \geq 2$. Let $N$ be the finite nerve of $W$, with
first barycentric subdivision $N'$, and let $K$ be the cone on $N'$.  For example, if $W$ is
generated by reflections in the sides of a right-angled hyperbolic $p$--gon $P$, then $N$ is a
circuit of $p$ edges, and $K$ is isomorphic to the barycentric subdivision of $P$.  For each $i
\in I$, let $G_i$ be a group of order $q_i$.  Each vertex of $K$ has a type $J$, where $J \subset
I$ is such that the group $W_J$ generated by $\{s_i \}_{i \in J}$ is finite.  For each $i$, let $K_i$ be the subcomplex of $K$ which is the closed star of the vertex of type $\{i\}$ in $N'$.  Let $G(K)$ be the
complex of groups where the vertex of $K$ with type $J$ has local group the direct product
\[\coprod_{i \in J} G_i\] and all monomorphisms are natural inclusions.  This complex of groups is
developable, with universal cover a right-angled building $X$ of type $(W,I)$.  The copies of $K$ in $X$ are called \emph{chambers}, and each copy of $K_i$ in $X$ is contained in $q_i$ distinct
chambers.  Moreover, the fundamental group of this complex of groups may be viewed as a uniform lattice
in $\Aut(X)$ (if all $q_i$ are finite).

Many other lattices for right-angled buildings (in any dimension) were obtained by promoting tree lattices, using
complexes of groups, in Thomas~\cite{Th2}.

\subsection{Kac--Moody buildings}\label{ss:KM_buildings}

Kac--Moody groups over finite fields $\F_q$ may be viewed as infinite-dimensional analogs of Lie groups.  See, for
example, Carbone--Garland~\cite{CG} and R\'emy--Ronan~\cite{RR}.  For any Kac--Moody group $\Lambda$ there are
associated (twin) buildings $X_+$ and $X_-$, constructed using twin $BN$--pairs $(B^+,N)$ and $(B^-,N)$ (see
Section~\ref{sss:classical} above).  The group $\Lambda$ acts diagonally on the product $X_+ \times X_-$, and for
$q$ large enough $\Lambda$ is a nonuniform lattice in $\Aut(X_+ \times X_-)$ (see~\cite{R1}).  A
\emph{Kac--Moody building} is a building which appears as one of the twin buildings for a Kac--Moody group.
Kac--Moody buildings are buildings, but unlike classical Euclidean buildings (see
Theorem~\ref{t:G_AutX}), non-isomorphic Kac--Moody groups may have the same building (R\'emy~\cite{R}).  One may
also study the \emph{complete Kac--Moody group} $G$, which is the closure of $\Lambda$ in the automorphism group
of one of its twin buildings.  Very few lattices in complete Kac--Moody groups are known.

\subsection{Davis--Moussong complexes}\label{ss:DM_complexes}

Given any Coxeter system $(W,S)$, the associated Davis--Moussong complex is a locally finite,
CAT(0), piecewise Euclidean polyhedral complex on which $W$ acts properly discontinuously and
cocompactly.  We describe a special case of this construction in dimension two.

Let $L$ be a connected, finite simplicial graph with all circuits of length at least $4$, and let
$k \geq 2$ be an integer.  The Coxeter system corresponding to this data has a generator $s_i$ of
order two for each vertex $v_i$ of $L$, and a relation $(s_i s_j)^{k} = 1$ if and only if the
vertices $v_i$ and $v_j$ are connected by an edge in $L$.  The Coxeter group defined by this
Coxeter system is denoted $W = W(k,L)$.  If $k = 2$, then $W$ is a right-angled Coxeter group.

For any such $W=W(k,L)$, Davis--Moussong constructed a CAT$(0)$ piecewise Euclidean complex $X =
X(2k,L)$ (see~\cite{Davis83,Moussong88}).  The cells of $X$ correspond to cosets in $W$ of
spherical subgroups of $W$, and in particular the $0$--cells of $X$ correspond to the elements of
$W$, viewed as cosets of the trivial subgroup.  Recall that a \emph{spherical subgroup} of $W$ is
a subgroup $W_T$ generated by some subset $T \leq S$, such that $W_T$ is finite.

The Davis--Moussong complex may be identified with (the first barycentric subdivision of) a
polygonal complex $X$ with all links $L$ and all $2$--cells regular Euclidean $2k$--gons.  The
group $W$ has a natural left action on $X$ which is properly discontinuous, cellular, and simply
transitive on the vertices of $X$.  Thus $W$ may be viewed as a uniform lattice in $\Aut(X)$. This
construction can also be carried out in higher dimensions, provided $L$ is a CAT$(1)$ spherical
simplicial complex.  In dimension $2$, where $L$ is a graph, this is equivalent to all circuits
having length at least $4$, by the Gromov Link Condition (Theorem~\ref{t:gromov_link} above).  If
$W$ is right-angled, then each apartment of a right-angled building of type $W$ is
isomorphic to the Davis--Moussong complex for $W$.

Davis--Moussong also found easy-to-verify conditions on $L$ such that $X(2k,L)$ may be equipped
with a CAT$(-1)$ piecewise hyperbolic structure.  In this way, some hyperbolic buildings (or
rather, their first barycentric subdivisions) may be constructed as Davis--Moussong complexes,
with the graph $L$ a one-dimensional spherical building.

\subsection{$(k,L)$--complexes}\label{ss:kL_complexes}

Let $L$ be a finite graph and $k$ an integer $\geq 3$.  A \emph{$(k,L)$--complex} is a polygonal complex $X$ such
that the link of each vertex of $X$ is $L$, and each $2$--cell of $X$ is a regular $k$--gon (usually but not
necessarily Euclidean).

\begin{figure}[ht]
\begin{center}
\includegraphics{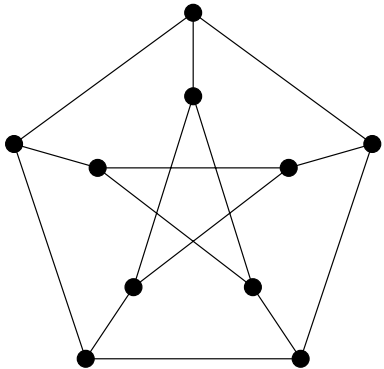}
\caption{Petersen graph}\label{f:petersen}
\end{center}
\end{figure}

Many polygonal complexes already described are $(k,L)$--complexes. For example, $2$--dimensional Euclidean or
hyperbolic buildings, with all links the same, are $(k,L)$--complexes where $L$ is a one-dimensional spherical
building.  The two-dimensional Davis--Moussong complexes described in
Section~\ref{ss:DM_complexes} above are barycentric subdivisions of $(k,L)$--complexes with $k \geq 4$ even.  An
example of a $(k,L)$--complex which is not a building or a Davis--Moussong complex is where $k$ is odd and $L$ is
the Petersen graph (Figure~\ref{f:petersen}).

There are simple conditions on the pair $(k,L)$ ensuring that a $(k,L)$--complex satisfies Gromov's Link Condition
(Theorem~\ref{t:gromov_link} above) and thus has nonpositive curvature.  Ballmann--Brin~\cite{BB1} showed that any
$(k,L)$--complex where $k$ and $L$ satisfy these conditions may be constructed in an inductive manner, by adding
$k$--gons to the previous stage without obstructions.  This construction is discussed in more detail in
Section~\ref{ss:localdata} below. Some $(k,L)$--complexes may also be constructed as universal covers of triangles of
groups, as done in~\cite{JLVV}.  In this case the fundamental group of the triangle of groups is a uniform
lattice. Constructions of uniform and nonuniform lattices as fundamental groups of complexes of groups are carried out for
certain highly symmetric $(k,L)$--complexes, including those with Petersen graph links, in Thomas~\cite{Th4}.

\subsection{CAT(0) cubical complexes}\label{ss:CAT0_cube}

Recall that a \emph{cubical complex} is a Euclidean polyhedral complex with  all $n$--cells isometric to the
Euclidean $n$--cube, and that a cubical complex $X$ is locally CAT(0) if and only if each vertex of $X$
is a flag simplicial complex (Theorem~\ref{t:flag} above).  Trees and products of trees are examples of
CAT(0) cubical complexes.

Groups of automorphisms of CAT$(0)$ cubical complexes are different in many ways from groups acting on the
Euclidean buildings discussed in Section~\ref{ss:Eucl_buildings} above.
Examples of discrete groups that act properly on CAT$(0)$ cube complexes include finitely generated Coxeter
groups~\cite{NibloReeves03}, many small cancellation groups~\cite{Wise04}, one-relator groups with torsion \cite{LauerWise}, many diagram groups, including Thompson's group
$F$ \cite{Farley03}, and groups acting properly on products of trees.

In this setting, the main geometric objects of study are hyperplanes, defined as follows. Consider two edges of a
CAT$(0)$ cube complex $X$ to be \emph{equivalent} if they are opposite edges of some $2$--cube. This relation generates
an equivalence relation whose equivalence classes are the \emph{combinatorial hyperplanes} of $X$. One can also
define \emph{geometric hyperplanes} of $X$ as unions of \emph{midplanes} of cubes, where a midplane of the cube
$C=[0,1]^n$ is a subset of the form \[ [0,1]\times\cdots\times[0,1] \times \{1/2\} \times [0,1] \times \cdots
\times[0,1]. \] Thus $C$ has $n$ midplanes, which intersect transversely at the barycenter of $C$. Given a
combinatorial hyperplane $H$, the corresponding geometric hyperplane is the union of all midplanes meeting the
barycenters of the edges of $H$.  Each geometric hyperplane is itself a CAT$(0)$ cubical complex, whose cubes are
midplanes of cubes of $X$. Each geometric hyperplane separates $X$ into two complementary components, called
\emph{halfspaces}. The properties of hyperplanes generalize the separation properties of edges in a tree. The main
new feature in higher dimensions, not present in trees, is that hyperplanes can have transverse intersections. In
fact, CAT$(0)$ cubical complexes have a rich combinatorial structure arising from the incidence and nesting
properties of hyperplanes.

Geometrically, the most significant subgroups in a group acting on a CAT(0) cubical complex are the
codimension--$1$ subgroups, which typically arise as stabilizers of hyperplanes. If a group $\Gamma$ has a finite
generating set $\mathcal{S}$, a subgroup $H \le \Gamma$ is \emph{codimension--$1$} provided that some neighborhood
of $H$ separates Cayley$(G,\mathcal{S})$ into at least two ``deep'' complementary components, where a component is
\emph{deep} if it contains elements arbitrarily far away from $H$. For instance, if $M$ is a $3$--manifold with an
immersed, incompressible surface $S$, then $\pi_1(S)$ is a codimension--$1$ subgroup of $\pi_1(M)$.

Sageev has shown (together with a result proved independently by Gerasimov and Niblo--Roller) that a finitely
generated group $\Gamma$ has a codimension--$1$ subgroup if and only if $\Gamma$ acts on a CAT$(0)$ cube complex
with no global fixed point \cite{Sageev95,Gerasimov97,NibloRoller98}. The cube complex produced by Sageev's
theorem is sometimes infinite dimensional and sometimes locally infinite.

Several representation-theoretic aspects of actions on trees extend naturally to actions on CAT$(0)$ cubical
complexes. If a topological group with Property (T) acts on a CAT$(0)$ cubical complex, then the action must have
a global fixed point. On the other hand, if a topological group $G$ acts metrically properly on a CAT$(0)$ cubical
complex $X$ then $G$ is a-T-menable \cite{NibloRoller98}. In particular, if $X$ is locally finite then any
discrete subgroup $\Gamma \le \Aut(X)$ is a-T-menable. Niblo--Reeves have also shown that if $X$ is any CAT$(0)$
cube complex, then every uniform lattice $\Gamma \le \Aut(X)$ is biautomatic \cite{NibloReeves98}.

\subsection{Systolic complexes}\label{ss:systolic}

Systolic complexes are a family of simplicial complexes whose geometry exhibits many aspects of nonpositive
curvature, yet which are not known to be CAT$(0)$.  A \emph{systolic complex} is a flag simplicial complex that is
connected and simply connected, such that the link of each vertex does not contain an isometric edge cycle of length
$4$ or $5$.  In~\cite{Chepoi2000}, Chepoi proved that a graph is the 1--skeleton of a systolic complex if and only if
it is a \emph{bridged graph}, which is a connected graph having no isometric edge cycles of length at least four.

Bridged graphs were introduced by Soltan--Chepoi~\cite{SoltanChepoi1983} and independently by
Farber--Jamison~\cite{FarberJamison1987} where they were shown to share certain convexity properties with CAT$(0)$ spaces. Their
geometric and algorithmic properties were studied for many years from the point of view of graph theory.  Systolic
complexes were rediscovered independently by Januszkiewicz--\Swiatkowski\ \cite{JS2}, Haglund \cite{Haglund2003}, and Wise, and
have subsequently been the subject of much study in geometric group theory.

A simplicial complex can be metrized in many ways, but the most natural metric, called the \emph{standard
piecewise Euclidean metric}, is given by declaring each simplex to be isometric to a regular Euclidean simplex
with all side lengths equal to $1$. In dimension $2$, a simplicial complex is systolic exactly when the its
standard piecewise Euclidean metric is CAT$(0)$. In higher dimensions being systolic is neither stronger nor
weaker than the standard piecewise Euclidean metric being CAT$(0)$.  A much more subtle question is whether a
systolic complex admits any piecewise Euclidean metric that is CAT$(0)$. No answer is known, but the answer is
generally expected to be negative.  

Systolic complexes do share many properties with CAT$(0)$ spaces. For example, any finite dimensional systolic
simplicial complex is contractible. As with CAT$(0)$ cubical complexes, any group acting properly discontinuously
and cocompactly on a systolic complex is biautomatic.  An interesting question is whether all systolic groups are in
fact CAT$(0)$ groups.

Systolic complexes are constructed in~\cite{JS2} as universal covers of simplices of groups, using the result that
a locally $6$--large complex of groups is developable.  The fundamental groups of these simplices of groups are
uniform lattices in the automorphism group of the universal cover.

\section{Properties of $\mathbf{X}$ and \textbf{Aut}$\mathbf(X)$}\label{s:properties_X}

The goal of this section is to understand the general structure of a polyhedral complex $X$ and
its full automorphism group $\Aut(X)$. For instance, how much local data is required in order to
uniquely determine $X$? What are the basic topological and group-theoretic properties of $\Aut(X)$?

\subsection{When do local data determine $X$?}\label{ss:localdata}

As seen in many examples in Section~\ref{s:examples} above, polyhedral complexes $X$ are often
constructed as universal covers of complexes of groups, and lattices in
$\Aut(X)$ are often fundamental groups of complexes of groups.  In each
case, the local structure of the universal cover is determined by the local structure of the
quotient space, together with the attached local groups of the complex of groups.
Thus it is critical to know how much local data is needed in order to uniquely specify a desired
polyhedral complex $X$.  To simplify matters, we focus on the special case when $X$ is a
$(k,L)$--complex (see Section~\ref{ss:kL_complexes} above).

\begin{ques} For a fixed $(k,L)$, is there a unique $(k,L)$--complex $X$?  If not, then what additional local data is
needed to determine $X$ uniquely? \end{ques}

If $L$ is the complete bipartite graph $K_{m,n}$, then in many cases there is a unique
$(k,L)$--complex $X$.  If $k=4$, this complex is the product of an $m$--valent and an $n$--valent
tree~\cite{Wise96}. If $k>4$ and either $k$ is even or $n=m$, the unique $(k,L)$--complex is
isomorphic to Bourdon's building $I_{p,q}$, a right-angled Fuchsian building, with $k = p$ and $L
= K_{q,q}$ (\cite{B1,Sw}; $I_{p,q}$ is discussed in Section~\ref{ss:hyp_buildings} above).  If
$k>4$ is odd and $n\ne m$ then there does not exist a $(k,L)$--complex.

On the other hand, when $L$ is the complete graph $K_n$ for $n \ge 4$, Ballmann--Brin~\cite{BB1}
and Haglund~\cite{H1} independently constructed uncountably many non-isometric
$(k,L)$--complexes.  We now
discuss these constructions.  As mentioned in Section~\ref{ss:kL_complexes}
above, simply connected nonpositively curved complexes can be constructed ``freely'' by building
successive balls outward from a given cell. Provided that certain obvious local obstructions do
not occur, we can glue in cells arbitrarily at each stage.  Ballmann--Brin showed that for many
choices of $k$ and $L$, every nonpositively curved $(k,L)$--complex can be constructed in this
manner~\cite{BB1}.

\begin{figure}
\begin{center}
\includegraphics{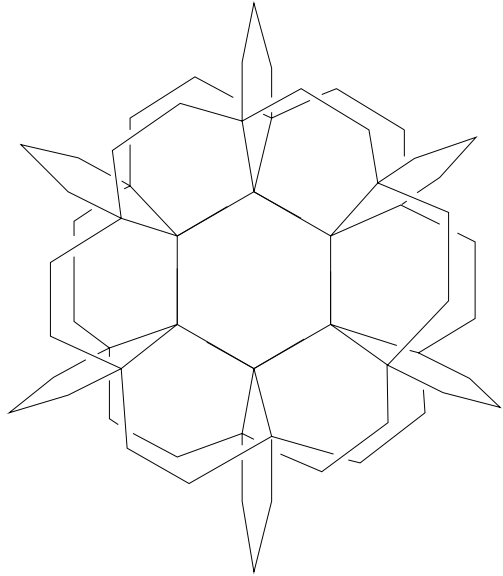}
\end{center}
\caption{If the union of the twelve hexagons
surrounding $A$ is a M\"obius band, then
the holonomy around the boundary of $A$ is nontrivial.}
\label{fig:BradyCrisp}
\end{figure}

In this inductive construction of a $(k,L)$--complex, choices may or may not arise.  Let us
consider the case when $k=6$ and $L=K_4$.  Then each $2$--cell of a $(k,L)$--complex $X$ is a
regular hexagon.  Each $1$--cell of $X$ is contained in three distinct hexagons.  Fix a $2$--cell $A$
of $X$, and consider the twelve surrounding $2$--cells which contain one of the six $1$--cells
bounding $A$.  These $2$--cells are arranged locally in two sheets, whose union is a band
surrounding $A$.  However, if one follows the sheets around the boundary of $A$, there are two
cases, depending on whether the union of the $12$ hexagons is an annulus, or is the M\"obius band
shown in Figure~\ref{fig:BradyCrisp}.

To describe and analyze this phenomenon, Haglund~\cite{H3} introduced the notion of \emph{holonomy}, which
measures the twisting of the $2$--cells neighboring a given $2$--cell $C$ as one traverses the
boundary cycle of $C$. In many cases, the choices of holonomies around each $2$--cell uniquely
determine the isomorphism type of a nonpositively curved $(k,L)$--complex. The existence of
holonomies depends on combinatorial properties of the graph $L$.

For instance, when $n \ge 4$, the complete graph $L=K_n$ admits nontrivial holonomies.  Roughly
speaking, Ballmann--Brin and Haglund constructed uncountably many $(k,K_n)$--complexes by showing
that, at each stage, a countable number of holonomies can be specified arbitrarily. In particular, $K_4$ has a
unique nontrivial holonomy, which is  illustrated in Figure~\ref{fig:BradyCrisp}. The unique
$(6,K_4)$--complex with trivial holonomies around every $2$--cell is the Cayley complex for the
presentation $\langle\,a,b \mid ba^2=ab^2 \,\rangle$, which defines the Geisking $3$--manifold
group. The unique $(6,K_4)$--complex with nontrivial holonomies around every $2$--cell is the
Cayley complex for the presentation $\langle\,a,b \mid aba^2=b^2 \,\rangle$, which is
$\delta$--hyperbolic (see \cite{BradyCrisp07} for more details).  On the other hand, the complete
bipartite graph $L=K_{m,n}$ admits only the trivial holonomy, which explains why there is a unique
$(k,L)$--complex in this case.

\Swiatkowski~\cite{Sw} considered $(k,L)$--complexes $X$ where $L$ is a trivalent graph and $X$ has \emph{Platonic symmetry}, that is, $\Aut(X)$ acts transitively on the set of flags (vertex, edge, face)
in $X$. He found elementary graph-theoretic conditions on $L$ that imply that
such an $X$ is unique.  Januszkiewicz--Leary--Valle--Vogeler~\cite{JLVV} classify Platonic
$(k,L)$--complexes $X$ in which $L$ is a complete graph.  Their main
results are for finite complexes $X$.

In general holonomies are not enough to uniquely determine a $(k,L)$--complex.
For instance, Haglund has observed that the Euclidean buildings for $\SL(3,\Q_p)$ and
$\SL \bigl(3,\F_p((t))\bigr)$ are $(3,L)$--complexes with the same link $L$ and the same holonomies.
Yet the buildings are not isomorphic, by Theorem~\ref{t:G_AutX} above.

Nonclassical buildings with given local structures have been studied
by Gaboriau--Paulin and Haglund--Paulin, who proved results analogous to
those for $(k,K_n)$--complexes and $(k,K_{m,n})$--complexes discussed
above.  If $q>4$ is a prime power,
Gaboriau--Paulin~\cite{GP} proved that for every hyperbolic Coxeter polygon $P$
with all vertex angles $\pi/6$, and for every prime power $q>4$, there exist uncountably many
hyperbolic buildings with chambers $P$ such that the links of vertices
are all isomorphic to the building for the projective plane over the
finite field $\F_q$.  On the other hand, Haglund--Paulin~\cite{HP2} showed that if $(W,I)$
is a right-angled Coxeter system and $(q_i)_{i\in I}$ is a collection
of cardinalities, then there exists a unique building $X$ of type
$(W,I)$ such that for each $i \in I$, each codimension one cell containing a vertex of type $\{i\}$ in $X$ is a face of $q_i$ distinct chambers.  This generalizes the result that Bourdon's building $I_{p,q}$ is the unique $(p,K_{q,q})$--complex.

In many cases it is still unknown how much local data is required to uniquely
specify a $(k,L)$--complex.






\subsection{Nondiscreteness of $\Aut(X)$}\label{ss:nondiscreteness}

Let $X$ be a locally finite, nonpositively curved polyhedral complex. The most basic question
about the locally compact group $G=\Aut(X)$ is whether or not it is discrete.    Recall that in
the compact-open topology, the group $G=\Aut(X)$ is nondiscrete exactly when, for each positive
integer $n$, there is an element $g_n \in G$, with $g_n$ fixing pointwise the ball of radius $n$
in $X$, and $g_n \not = {\rm Id}$.  The theory of lattices in a discrete group is trivial, hence this
issue is of crucial importance.  We again focus on the case of $(k,L)$--complexes (see
Section~\ref{ss:kL_complexes} above).

\begin{ques}
Given a $(k,L)$--complex $X$, is $G=\Aut(X)$ discrete?
\end{ques}

The answer is known in certain cases, and is closely related to the notion of a flexible complex.

\begin{dfinition} A complex $X$ is \emph{flexible} if there exists $\phi\in \Aut(X)$ such that $\phi$ fixes the
star of some vertex in $X$ but $\phi\neq\Id$.\end{dfinition} \medskip

Flexibility was introduced by Ballmann--Brin in \cite{BB1}. If $X$ is locally finite and not flexible, then the
stabilizer of any vertex $v \in X$ is finite, since an automorphism of $X$ that fixes $v$ is uniquely determined
by its action on the link of $v$. In particular, $\Aut(X)$ is discrete if $X$ is not flexible. The following
result is nearly immediate from the definition of flexibility.

\begin{theorem}[Discreteness criterion]\label{t:discreteness} If the graph $L$ is not flexible, then
no $(k,L)$--complex $X$ is flexible, and $\Aut(X)$ is discrete. \end{theorem}

Theorem~\ref{t:discreteness} has the following converse when $X=X(2k,L)$ is the Davis--Moussong complex for the
Coxeter group $W=W(k,L)$, discussed in Section~\ref{ss:DM_complexes} above. The result was proved independently by
Haglund and \Swiatkowski\ in the case that $X$ is $2$--dimensonal \cite{H2,Sw}, and was  extended to arbitrary
Coxeter systems by Haglund--Paulin \cite{HP1}.

\begin{theorem}[Nondiscreteness criterion]\label{t:nondiscreteness} Suppose $L$ is a finite simplicial graph and
$k\ge 2$. Let $X=X(2k,L)$ be the Davis--Moussong complex for the Coxeter group $W=W(k,L)$. If $L$ is flexible then
$\Aut(X)$ is nondiscrete. \end{theorem}

The proof of Theorem~\ref{t:nondiscreteness} relies on the fact that Davis--Moussong complexes have numerous
symmetries.  For other $(k,L)$--complexes, particularly those with $k$ odd, much less is known. It is not clear
whether this reflects the limitations of our techniques, or actual differences in behavior for $k$ odd and $k$
even.

\subsection{Simplicity and nonlinearity}\label{ss:simplicity_nonlinearity}

Let $X$ be a locally finite, nonpositively curved polyhedral complex, with locally compact automorphism group
$G=\Aut(X)$.  In this section we discuss whether two basic group-theoretic properties, simplicity
and (non)linearity, hold for $G$.  We assume for this section that $G=\Aut(X)$ is nondiscrete.

\begin{ques} When is $G=\Aut(X)$ a simple group? \end{ques}

For $X$ a locally finite regular or biregular tree, Tits~\cite{T1} proved simplicity of the group $\Aut_0(X)$ of
type-preserving automorphisms of $X$ (which is finite index in the full automorphism group $G=\Aut(X)$).
Haglund--Paulin~\cite{HP1} showed that various type-preserving automorphism groups in several higher-dimensional
cases are simple.  We note that the method of proof of these results lies in geometric group theory.

We say that a group $G$ is {\em linear} if it has a faithful representation $G\to \GL(n,K)$ for some field $K$.
On the question of linearity, suppose $X$ is a classical Euclidean building, associated to the algebraic group
$G$ over a local nonarchimedean field $K$ (see Section~\ref{sss:classical} above).   Theorem~\ref{t:G_AutX}
above says that if $\ch(K) = 0$ then $G$ is finite index in $\Aut(X)$ (and if $\ch(K) = p > 0$ then $G$ is
cocompact in $\Aut(X)$).  By inducing, we see in particular that when $\ch(K) = 0$, the group $\Aut(X)$
has a faithful linear representation over $K$.  On the other hand, for several higher-dimensional complexes $X$
which are not classical buildings, Haglund--Paulin~\cite{HP1} proved that the full automorphism group $\Aut(X)$
has no such faithful linear representation.  For $\dim(X)=2$, we pose the following question:

\begin{prob} Find conditions on the link $L$ so that a $(k,L)$--complex $X$ has
linear automorphism group.  \end{prob}

\noindent Haglund~\cite{H4} has recently shown that $\Aut(X)$ is nonlinear for certain Fuchsian buildings $X$ (see
Section~\ref{ss:hyp_buildings} above). Is it possible that linearity of $\Aut(X)$ characterizes those $X$ which
are classical Euclidean buildings, among all nonpositively curved $X$?

\section{Comparisons with linear groups}\label{s:comparisons}

While one expects some of the phenomena and results from the theory of linear
groups $G\subset \GL(n,\C)$ to hold for the group $G=\Aut(X)$ and its lattices,
most of the methods from that theory are unavailable in this new context.  There
are no eigenvalues or traces.  There are no vectors to act on.  It therefore
seems important to attack such questions, as they will (hopefully) force us to
come up with new methods and tools.

One new approach to the study of automorphism groups of nonpositively curved
polyhedral complexes is the structure theory of totally disconnected locally
compact groups (see the survey~\cite{Willis}).  An example
of this approach is the computation of the flat rank of automorphism groups of
buildings with sufficiently transitive actions~\cite{BRW}.

\subsection{Some linear-type properties}\label{ss:linear_type}

One of the basic properties of linear groups $G$ is the {\em Tits alternative}: any finitely generated
linear group either contains a nonabelian free group or has a solvable subgroup of finite index (see~\cite{T2}).
The following problem is well known.

\begin{prob} \label{prob:tits} Let $X$ be a nonpositively curved polyhedral complex. Prove that finitely generated
subgroups of $G=\Aut(X)$ satisfy the Tits alternative. \end{prob}

When $X$ is a CAT$(-1)$ space, uniform lattices in $G=\Aut(X)$ are word-hyperbolic, and thus satisfy the Tits
alternative (Gromov~\cite{Gr}).  The usual ping-pong argument for the Tits alternative requires strong
expanding/contracting behavior for the action of isometries of $X$ on the visual boundary $\partial X$.  The
difficulty with Problem \ref{prob:tits} lies in the fact that if $X$ is just nonpositively curved, rather than
negatively curved, this behavior on $\partial X$ is not strong enough to immediately allow for the
usual ping-pong argument to work.

The Iwasawa decomposition $KAN$ of a semisimple Lie group $G$ plays a fundamental role in the representation
theory of $G$.  Here, $K$ is a compact subgroup, $A$ is abelian and $N$ is nilpotent.  In the topology on
$G=\Aut(X)$, where $X$ is a locally finite polyhedral complex, the stabilizers of vertices are maximal compact subgroups.

\begin{ques} For which $X$ does $G=\Aut(X)$ have a $KAN$ structure?
\end{ques}

Answering this question might be a first step towards investigating various analytic properties of $X$, the group $G=\Aut(X)$, and its lattices.  For instance, random walks on classical buildings have been studied using the representation theory of the associated algebraic group (see, for example, Cartwright--Woess~\cite{CW} and Parkinson~\cite{Pa}), but for more general complexes $X$ this machinery is not available.

Kazhdan proved that simple Lie groups $G$ have {\em property \textup{(}T\textup{)}}: the trivial representation is isolated in the
unitary dual of $G$ (see, for example,~\cite{M}). Ballmann--\Swiatkowski~\cite{BS}, \.Zuk~\cite{Z}, and
Dymara--Januszkiewicz~\cite{DJ} have proven that many $G=\Aut(X)$ satisfy this important property.

\begin{ques} For which $X$ does $\Aut(X)$ have Property (T)?
\end{ques}

\noindent We remark that a locally compact topological group $G$ has property (T) if and only if any of its
lattices has property (T).

One of the deepest theorems about irreducible lattices $\Gamma$ in higher rank semisimple Lie groups is Margulis's
Normal Subgroup Theorem (see~\cite{M}), which states that any normal subgroup of $\Gamma$ is finite or has finite
index in $\Gamma$.

\begin{ques}\label{ques:normalsubgp} For which $X$ does a normal subgroup theorem hold for $\Aut(X)$?\end{ques}

\noindent Such a theorem has been shown for products of trees by Burger--Mozes~\cite{BM}.

Recall that the {\em Frattini subgroup} $\Phi(\Gamma)$ of a group $\Gamma$ is the intersection of all maximal
subgroups of $\Gamma$.  Platonov~\cite{Pl} proved that $\Phi(\Gamma)$ is nilpotent for every finitely generated
linear group.  Ivanov~\cite{I} proved a similar result for mapping class groups.  I. Kapovich~\cite{K} proved that
$\Phi(\Gamma)$ is finite for finitely generated
subgroups of finitely generated word-hyperbolic groups.

\begin{prob}
Compute the Frattini subgroup $\Phi(\Gamma)$ for finitely generated subgroups $\Gamma<\Aut(X)$.
\end{prob}

Part of the fascination of lattices in $\Aut(X)$ is that they exhibit a mixture of rank one and higher rank
behavior.  Ballmann--Eberlein (see~\cite{E}) defined an invariant $\rank(\Gamma)$, called the \emph{rank} of
$\Gamma$, which is defined for any finitely generated group $\Gamma$ as follows.   Let $\Gamma_i$ denote the set
of elements $g\in\Gamma$ so that the centralizer of $g$ contains $\Z^d$ for some $d\leq i$ as a finite index
subgroup.  Let $r(\Gamma)$ be defined to be the smallest $i$ so that $\Gamma$ is a finite union of translates
$$\Gamma=g_1\Gamma_i \cup \cdots \cup g_n\Gamma_i$$ for some $g_j\in\Gamma$.  Then define $\rank(\Gamma)$ to be
the maximum of $r(\Gamma')$, where $\Gamma'$ runs over all finite index subgroups of $\Gamma$.

Work of Prasad--Raghunathan shows that this notion of rank agrees with the classical one for arithmetic lattices.
Ballmann--Eberlein~\cite{BE} proved that the rank of the fundamental group of a complete, finite volume,
nonpositively curved manifold $M$ equals the geometric rank of the universal cover of $M$.  Since centralizers of
infinite order elements in word-hyperbolic groups $\Gamma$ are virtually cyclic, it is clear that
$\rank(\Gamma)=1$ in these cases. Thus for nonpositively curved, connected, simply-connected $2$--complexes $X$,
lattices in $\Aut(X)$ can have rank one and also rank two (the latter, for example, when $X$ is a classical
Euclidean building, discussed in Section~\ref{sss:classical} above).

 \begin{prob}
 Compute $\rank(\Gamma)$ for lattices $\Gamma<\Aut(X)$.
 \end{prob}

A basic property of any finitely generated linear group is that it is residually
finite. In contrast, there are lattices $\Gamma$ in $G=\Aut(X)$ that are not
residually finite.  Indeed, Burger--Mozes \cite{BM} have constructed, in the
case when $X$ is a product of simplicial trees, lattices which are simple
groups.  Wise had earlier constructed lattices for such $X$ that are not
residually finite~\cite{W1}.  Kac--Moody lattices are also simple, and their buildings have
arbitrarily large dimension (see~\cite{CapraceRemy}).

\begin{prob} Construct a lattice $\Gamma$ in $G=\Aut(X)$ which is a simple
group, and where $X$ is not a product of trees. \end{prob}

For residual finiteness, a key case is Bourdon's building $I_{p,q}$ (see Section~\ref{ss:hyp_buildings} above)
whose $2$--cells are right-angled hyperbolic $p$--gons. Wise~\cite{W} has shown that fundamental groups of
polygons of finite groups, where the polygon has at least 6 sides, are residually finite.  Thus there are
residually finite uniform lattices for $I_{p,q}$, $p \geq 6$, but the question of whether
every uniform lattice in $I_{p,q}$ is residually finite is completely
open for $p=5$, that is, for pentagons.  The question of residual finiteness of uniform lattices is open
even for triangular hyperbolic buildings (see~\cite{KV}).

\begin{ques}
Which lattices $\Gamma<G=\Aut(X)$ are residually finite?
\end{ques}

A related but broader problem is as follows.  Most of the known CAT$(0)$ groups are residually finite, hence virtually torsion-free.  As remarked in
Section~\ref{ss:cxs_of_groups} above, to date most applications of the theory of complexes of groups have used only
simple complexes of groups.  Now, if the fundamental group $\G$ of a complex of groups $G(Y)$ is virtually
torsion-free, then $G(Y)$ has a finite cover $G(Y') \to G(Y)$ where all local groups of $G(Y')$ are trivial, hence
$G(Y')$ is a simple complex of groups.  

\begin{prob}[Haglund] Find a nonpositively curved complex of groups $G(Y)$ which is not finitely
covered by a simple complex of groups.  Do this in the negatively curved setting as well.  Is there a CAT$(0)$ group
$\G$ which is not virtually the fundamental group of \emph{any} (nonpositively curved) simple complex of groups?
\end{prob}

\subsection{Rigidity}\label{ss:rigidity}

Automorphism groups $G$ of nonpositively curved polyhedral complexes $X$, and lattices $\Gamma<G$, are natural
places in which to study various rigidity phenomena, extending what we know in the classical, algebraic cases.  A
first basic problem is to prove strong (Mostow) rigidity. In other words, one wants to understand the extent to
which a lattice $\Gamma$ in $G$ determines $G$.

\begin{prob}[Strong rigidity] Let $X_1$ and $X_2$ be nonpositively curved polyhedral complexes, and let $\Gamma_i$
be a lattice in $G_i=\Aut(X_i), i=1,2$.  Find conditions on the $X_i$ which guarantee that any abstract group isomorphism $\phi\colon\Gamma_1\to \Gamma_2$ extends to an isomorphism $G_1\to G_2$. Further, determine when any two
copies of\/ $\Gamma_i$ in $G_i$ are conjugate in $G_i$. \end{prob}

\noindent Some assumptions on the $X_i$, for example that every $1$--cell is contained in a $2$--cell, are needed to
rule out obvious counterexamples.

A harder, more general problem is to prove quasi-isometric rigidity.

\begin{prob}[Quasi-isometric rigidity] Compute the quasi-isometry groups of nonpositively curved polyhedral
complexes $X$.  Prove {\em quasi-isometric rigidity} theorems for these complexes; that is, find conditions on $X$
for which: \begin{enumerate} \item Any quasi-isometry of $X$ is a bounded distance from an isometry
(automorphism), and \item Any finitely-generated group quasi-isometric to $X$ is (a finite extension of) a cocompact lattice in $\Aut(X)$.
\end{enumerate} \end{prob}

A standard trick due to Cannon--Cooper shows that (1) implies (2). It is also immediate from Mostow's original
argument that (1) implies strong rigidity.  Quasi-isometric rigidity was proven in the case of Euclidean buildings
by Kleiner--Leeb \cite{KL}. Bourdon--Pajot~\cite{BP} proved quasi-isometric rigidity for Bourdon's building
$I_{p,q}$, and Xie~\cite{X} generalized this to Fuchsian buildings (see Section~\ref{ss:hyp_buildings}).  One would expect that higher-dimensional
buildings would be more rigid, and indeed they seem to be harder to construct, so they might be a good place to
look for rigidity phenomena.

Another kind of rigidity problem is the following:

\begin{prob} Suppose $X_1$ and $X_2$ are locally finite, connected, simply-connected 2-complexes, such that for $i
=1,2$, the group $\Aut(X_i)$ acts cocompactly on $X_i$. If $\Aut(X_1)$ is isomorphic to $\Aut(X_2)$, is $X_1$
isometric to $X_2$?\end{prob}

A variety of other rigidity phenomena from Riemannian geometry have natural analogs in this context.  Examples
include rank rigidity, hyperbolic rank rigidity, minimal entropy rigidity, and marked length spectrum rigidity.  A
rank rigidity theorem for nonpositively $2$--complexes was proven by Ballmann--Brin in \cite{BB2}.





\subsection{Geometry of the word metric}

One of the few results about the geometry of the word metric for nonuniform lattices is the theorem of
Lubotzky--Mozes--Raghunathan \cite{LMR}, which we now discuss.

Let $G$ be a semisimple Lie group over $\R$ (respectively, over a nonarchimedean local field $K$), and let $X$ be
the associated symmetric space (respectively, Euclidean building).  Thus $X$ is a nonpositively curved Riemannian
manifold (respectively, simplicial complex) on which $G$ acts by isometries.  Let $\Gamma$ be a lattice in $G$. If $K$ is
nonarchimedean and $G$ has rank one over $K$ then nonuniform lattices in $G$ are not finitely generated. On the
other hand, when $G$ is either a real Lie group or a nonarchimedean group with $K$--rank at least $2$,  then all
lattices $\Gamma$ in $G$ are finitely generated. In this section, we consider only the finitely-generated case,
and we endow $\Gamma$ with the word metric for a finite generating set.

If $\Gamma$ is uniform, then the natural map $\psi\colon \Gamma \to X$ sending $\Gamma$ to any of its orbits is a
quasi-isometry. When $\Gamma$ is nonuniform, the orbit map is never a quasi-isometry, since the quotient
$\Gamma\backslash X$ is noncompact. When $G$ has real rank one, the map $\psi$ is not even a quasi-isometric
embedding, as can be seen by considering any nonuniform lattice acting on real hyperbolic space.  In this case,
the maximal parabolic subgroups of $\Gamma$ are exponentially distorted in $X$.

The theorem of Lubotzky--Mozes--Raghunathan \cite{LMR} states that, when $G$ has real rank (respectively, $K$--rank) at
least $2$, then $\psi$ is indeed a quasi-isometric embedding.  Each of the known proofs of this result is heavily
algebraic, depending on the structure of matrix groups.  Thus the following problem presents an interesting
challenge, even in terms of giving a geometric proof in the (nonarchimedean) algebraic case.

\begin{ques} Let $\Gamma$ be a finitely generated nonuniform lattice in the automorphism group of a nonpositively curved polyhedral
complex $X$.  When is the natural map $\psi\colon \Gamma\to X$, sending $\Gamma$ to
any of its orbits, a quasi-isometric embedding? \end{ques}

When $X$ is a product of trees, $\psi$ need not be a quasi-isometric embedding.  When $X$ is not a product of
trees, is $\psi$ always a quasi-isometric embedding?

\subsection{Dynamics}

Let $G$ be (any) locally compact topological group, equipped with Haar measure, and let $\Gamma$ be a lattice in
$G$.  Then $G$ acts on the left on $G/\Gamma$, preserving the finite measure on $G/\Gamma$
induced by the Haar measure on $G$.  We thus obtain an action of every closed subgroup $H<G$ on $G/\Gamma$.  It is
a basic question understand these dynamical systems, in particular to determine when the action of $H$ on
$G/\Gamma$ is {\em ergodic}; that is, when every $H$--invariant set has zero or full measure.  When $G$ is a
semisimple Lie group with no compact factors, and $\Gamma$ is an irreducible lattice in $G$, {\em Moore's
Ergodicity Theorem} (see~\cite{Zim}) states that the $H$--action on $G/\Gamma$ is ergodic if and only if $H$ is
noncompact.

Now let $X$ be a simply connected, locally finite polyhedral complex of nonpositive curvature.  Equip
$G=\Aut(X)$ with left-invariant Haar measure, and let $\Gamma$ be a lattice in $G$.

\begin{prob}
\label{prob:Moore}
Determine which closed subgroups of $G=\Aut(X)$ act ergodically on $G/\Gamma$.
\end{prob}

One reason we consider Problem \ref{prob:Moore} to be worthwhile is that the usual method of proving Moore's
Ergodicity Theorem uses the unitary representation theory of $G$.  We thus believe that, apart from being
interesting in its own right, attempts to solve Problem \ref{prob:Moore} will require us either to find new
approaches to Moore's theorem, or to develop the unitary representation theory of $G=\Aut(X)$.

\section{Lattices in \textbf{Aut}$\mathbf{(X)}$}\label{s:lattices}

In this section we consider properties of the lattices in $\Aut(X)$ themselves. Some lattice properties have
already been mentioned in Section~\ref{s:comparisons} above, on comparisons with linear groups.  Here, we discuss topics where new phenomena, contrasting with classical cases, have already been observed, and where the known techniques of proof are combinatorial or geometric in flavor.

\subsection{Existence and classification theorems}\label{ss:existence_classification}

Given a locally compact group $G$, the most basic question in the lattice theory of $G$ is whether $G$ admits a
uniform or nonuniform lattice.

For algebraic groups, the existence of both uniform and nonuniform lattices was settled by Borel and others, using
arithmetic constructions (see the final paragraph of Section~\ref{sss:classical} above).  For automorphism groups of
trees, precise conditions are known for the existence of both uniform lattices (Bass--Kulkarni~\cite{BK}) and
nonuniform lattices (Bass--Carbone--Rosenberg, in~\cite{BL1}).   In Section~\ref{s:examples} above, for each example
$X$ of a polyhedral complex, we described known constructions of lattices in $G=\Aut(X)$.  These constructions are
non-arithmetic, for $X$ not a classical building.  The following question is still largely open.

\begin{ques} When does $G=\Aut(X)$ admit a uniform lattice? A nonuniform lattice? \end{ques}

\noindent A special case of this question is:

\begin{ques} For which positive integers $k \geq 3$ and finite simplicial graphs $L$ does the automorphism group of
a $(k,L)$--complex $X$  admit lattices?\end{ques}

Once one establishes the existence of lattices in a given $G=\Aut(X)$, the next problem is to classify all such
lattices.  We discuss commensurability of lattices in Section~\ref{ss:commens} below.  An even more fundamental
question is:

\begin{prob}
\label{prob:classify:lattices}
Classify lattices in $G=\Aut(X)$ up to conjugacy.
\end{prob}

We note that in the case of real Lie groups, classification theorems are difficult. For $\SO(3,1)$, for example, the classification is precisely the classification of all finite volume, complete hyperbolic orbifolds.
On the other hand, for higher rank real (and $p$-adic) semisimple Lie groups, Margulis's arithmeticity theorem (see~\cite{M}) states that all lattices are arithmetic, and arithmetic lattices can in some sense be classified
(although this is also not easy).  So, even solving Problem \ref{prob:classify:lattices} in any special case, for example for
specific hyperbolic buildings, would be of great interest.

\subsection{Commensurability}\label{ss:commens}

One of the basic problems about a locally compact topological group $G$ is to
classify its lattices up to commensurability. Recall that two lattices
$\Gamma_1,\Gamma_2 \le G$ are \emph{commensurable in $G$} if there exists $g \in
G$ so that $g \Gamma_1 g^{-1} \cap \Gamma_2$ has finite index in both $g
\Gamma_1 g^{-1}$ and $\Gamma_2$. Since covolume is multiplicative in index, two
commensurable lattices have covolumes that are commensurable real numbers, that is, they have a rational ratio.

\begin{prob}\label{prob:commens} Classify lattices in $G=\Aut(X)$ up to
commensurability.  As a subproblem, find commensurability invariants of lattices.
\end{prob}

If $G$ is an algebraic group of rank at least two over a nonarchimedean local field $K$, then there exist
noncommensurable arithmetic lattices in $G$. If $G$ is a rank one simple real Lie group, then lattices are again not
all commensurable, as there exist both arithmetic and nonarithmetic lattices.

For $G=\Aut(X)$, commensurability of uniform lattices is strikingly different. When $X$ is a locally finite tree,
Leighton proved in \cite{Leighton82} that all torsion-free uniform lattices in $G=\Aut(X)$ are commensurable. The
torsion-free hypothesis was removed by Bass--Kulkarni in \cite{BK}, establishing that there is at most one
commensurability class of uniform lattices in the tree case. Haglund~\cite{H4} has shown the same result for many
Fuchsian buildings (see Section~\ref{ss:hyp_buildings}). He has also found a sufficient condition for a uniform
lattice in the automorphism group of a Davis--Moussong complex $X$ (see Section~\ref{ss:DM_complexes}) to be
commensurable to the corresponding Coxeter group $W$.  As specific instances of Problem~\ref{prob:commens}, we
have:

\begin{prob} Suppose $X$ is a $(k,L)$--complex.  Find conditions on $L$ such that all uniform lattices in $\Aut(X)$
are commensurable, and find examples of such $L$.
\end{prob}

\noindent and on the other hand:

\begin{prob}[Haglund] Find a Gromov-hyperbolic CAT$(0)$ complex $X$ such that $\Aut(X)$ admits two non-commensurable uniform
lattices. \end{prob}

For nonuniform lattices in $G=\Aut(X)$, the situation
seems much wilder. Even in the tree case, there seems to be a great
deal of flexibility in the construction of nonuniform lattices. For
instance, Farb--Hruska \cite{FH} have shown that when $X$ is the biregular tree
there are uncountably many commensurability classes of nonuniform
lattices in $G=\Aut(X)$ with any given covolume $v>0$ . To
prove this result, they construct several new commensurability
invariants, and then evaluate them on lattices constructed using
graphs of groups.

A similar result holds when $X$ is a right-angled building (see Section~\ref{ss:rab}), by work of Thomas \cite{Th2}.  Lattices in right-angled hyperbolic buildings, such as Bourdon's building $I_{p,q}$, are known to exhibit higher-rank phenomena, such as quasi-isometric rigidity
(see \cite{BP} and Section~\ref{ss:rigidity} above).  In contrast, Thomas' theorem indicates a similarity of
these lattices with tree lattices.  In fact, Thomas proves this theorem by constructing a functor
that takes tree lattices to lattices in right-angled buildings.  This functor preserves many features of the lattice.

The most important commensurability invariant of a group $\Gamma$ inside a group $G$
is the {\em commensurator} $\Comm_G(\Gamma)$ of $\Gamma$ in $G$, defined by
$$\Comm_G(\Gamma):=\{\,g\in G\mid \Gamma\cap g\Gamma g^{-1}\ \mbox{has finite index in both $\Gamma$ and $g\Gamma g^{-1}$}\,\}.$$
Margulis proved that a lattice $\Gamma$ in a semisimple Lie group $G$ is arithmetic if and only if
$\Comm_G(\Gamma)$ is dense in $G$ (see~\cite{Zim}).  Lubotzky proposed this density property as a definition of ``arithmeticity" when $G=\Aut(X)$.

\begin{prob} For lattices $\Gamma$ in $G=\Aut(X)$, compute $\Comm_G(\Gamma)$.  Determine
whether or not $\Comm_G(\Gamma)$ is dense in $G$.
\end{prob}

When $X$ is a tree, density of commensurators of uniform lattices was proved by Bass--Kulkarni \cite{BK}
and Liu \cite{L}.  Haglund established density of commensurators of uniform lattices for many Davis--Moussong
complexes in\cite{H2}, and Haglund~\cite{H5} and independently Barnhill--Thomas~\cite{BT} have recently shown the same result for right-angled buildings.  For commensurators of nonuniform lattices, however, even for trees very little is known (see~\cite{BL1}).

\subsection{Finiteness properties of lattices}

Uniform lattices in $G=\Aut(X)$ are always finitely generated, for obvious reasons.
However, nonuniform lattices need not be finitely generated.

\begin{ques}\label{q:generation}
For which $G=\Aut(X)$ are all nonuniform lattices non--finitely generated?  Do there exist
$G$ which admit both finitely generated and non--finitely generated nonuniform
lattices?
\end{ques}

Higher rank algebraic groups, such as $G=\SL\bigl(3,( \F_q((t)) \bigr)$, have Kazhdan's Property (T) (see Section~\ref{ss:linear_type} above).  Furthermore, Property (T) is inherited by lattices,
and all countable groups with Property (T) are finitely generated.
Therefore lattices in higher rank groups are all finitely generated.

On the other hand, if $X$ is a tree, every nonuniform lattice in $\Aut(X)$
is non--finitely generated \cite{BL1}.  Thomas' functor
mentioned in Section~\ref{ss:commens} above implies that many nonuniform lattices
in right-angled hyperbolic buildings are non--finitely generated as well.

\begin{conjecture}
\label{conj:fg} Let $\Gamma$ be a nonuniform lattice in $G=\Aut(X)$,
where $X$ is any right-angled hyperbolic building.  Then $\Gamma$ is
not finitely generated.
\end{conjecture}

We are starting to believe that finite generation of nonuniform
lattices in $2$--complexes is actually a miracle, and could even
characterize the remarkable nonuniform lattices in algebraic groups
in characteristic $p>0$.   Even these lattices are not finitely presentable, and so we
make the following:

\begin{conjecture}
If $\Gamma$ is a nonuniform lattice in $G=\Aut(X)$, where $X$ is a locally finite polyhedral complex, then $\Gamma$ is not
finitely presentable.
\end{conjecture}

\subsection{Covolumes}

\label{s:covolumes}

One of the more striking ways in which the study of lattices in $\Aut(X)$ diverges from the case of lattices in semisimple Lie groups is the study of covolumes of lattices in a fixed $\Aut(X)$. New phenomena are seen to occur right away, and much remains to be understood.

\begin{prob}\label{p:covolumes} Given $G=\Aut(X)$ with Haar measure $\mu$, describe the
set of covolumes
\[
   \mathcal{V}(G) := \bigl\{\, \mu(\G \backslash G) \bigm|
   \mbox{$\G$ is a lattice in $G$}\,\bigr\}
\]
\end{prob}

\noindent Note that $\mathcal{V}(G)$ is a set of positive real
numbers.

If $G$ is a non-compact simple real Lie group, such as $\PSL(n,\R)$,
then the set $\mathcal{V}(G)$ has positive lower bound
(Kazhdan--Margulis,~\cite{KM}) and in most cases is discrete
(see~\cite{Lu} and the references therein).  If $G$ is a higher-rank
algebraic group over a nonarchimedean local field, such as
$\PSL(n,\Q_p)$ with $n \geq 3$, the strong finiteness result of
Borel--Prasad~\cite{BPr} implies that for any $c > 0$, there are only finitely
many lattices in $G$ with covolume less than~$c$. Hence
$\mathcal{V}(G)$ is discrete, has positive lower bound, and for any
$v \in \mathcal{V}(G)$ there are only finitely many lattices of
covolume~$v$.

The set of covolumes for tree lattices is very different.  Suppose
$G$ is the group of automorphisms of a regular locally finite
tree.  Then, for example, Bass--Kulkarni~\cite{BK} showed that
$\mathcal{V}(G)$ contains arbitrarily small elements, by
constructing a tower of uniform lattices (see Section~\ref{s:towers}
below).  Bass--Lubotzky~\cite{BL1} showed that the set of nonuniform
covolumes is $(0,\infty)$.

A few higher-dimensional nonclassical cases have been studied.
In~\cite{Th2} and~\cite{Th3}, Thomas considered covolumes for,
respectively, right-angled buildings and certain Fuchsian buildings (see Sections~\ref{ss:rab} and~\ref{ss:hyp_buildings} respectively).  In both these settings, $\mathcal{V}(G)$ shares properties, such as
nondiscreteness, with covolumes of tree lattices, even though such
buildings also have some rigidity properties typical of classical cases
(see Section~\ref{ss:rigidity} above).  Little is known about covolumes for $X$ not a building.
In~\cite{Th4}, the class of $(k,L)$--complexes $X$ of Platonic symmetry (introduced by
\Swiatkowski~\cite{Sw}; see Section~\ref{ss:kL_complexes}) is considered. A sample result is that
if $k \geq 4$ is even, and $L$ is the Petersen graph, then
$\mathcal{V}(G)$ is nondiscrete. Many cases are completely open.

From a different point of view, Prasad~\cite{P} gave a computable
formula for the covolumes of lattices $\G$ in algebraic groups $G$ over nonarchimedean local fields.  This formula is in terms of discriminants of field extensions and numbers of roots. If $\G$ is
viewed instead as a lattice in $\Aut(X)$, where $X$ is the building
associated to the algebraic group $G$, we have also Serre's more
geometrically-flavored formula for the covolume of $\G$, stated in
the introduction.

\begin{ques} Can Serre's geometric formula for covolumes tell us anything new about
lattices in classical cases? \end{ques}

\noindent More generally, using Serre's geometric formula, in~\cite{Th1}
Thomas established a computable number-theoretic restriction on the
set of covolumes of uniform lattices, for all locally finite $X$
with $G=\Aut(X)$ acting cocompactly, in all dimensions.

\begin{prob}\label{p:unif_covolumes} Suppose $v>0$ satisfies the
restriction of~\cite{Th1}. Construct a uniform lattice in $G$ of
covolume $v$, or show that such a lattice does not exist. Also, find
the cardinality of the set of uniform lattices of covolume $v$.  For nonuniform lattices, the same questions for any $v>0$.
\end{prob}

\noindent This problem was solved for right-angled buildings (see Section~\ref{ss:rab})
in~\cite{Th2}.

The properties of the set of volumes of hyperbolic three-manifolds are well-understood (see~\cite{Th}), and one could investigate whether similar properties hold for volumes of lattices in $\Aut(X)$.  For instance, for every nonuniform lattice $\Gamma$ in $\SO(3,1)$, there is a sequence of uniform lattices with covolumes converging to that of $\Gamma$, obtained by Dehn surgery.  This gives a surjective homomorphism from $\Gamma$ to each of these uniform lattices.  It is not known whether any nonuniform lattices in $\Aut(X)$ surject onto uniform lattices.

\subsection{Towers}\label{s:towers}

The study of towers of lattices in closely related to covolumes (Section~\ref{s:covolumes} above).  A \emph{tower} of lattices in a locally compact group $G$ is an infinite strictly ascending
sequence
\[\G_1 < \G_2 < \cdots < \G_n < \cdots < G\]
where each $\G_n$ is a lattice in $G$.

\begin{ques}\label{q:towers} Does $G=\Aut(X)$ admit a tower of (uniform or nonuniform) lattices?\end{ques}

If $G$ admits a tower, then the covolumes of lattices in this tower tend to zero, hence the set $\mathcal{V}(G)$ of covolumes does not have positive lower bound.  It follows that in classical (algebraic) cases, $G$ does not admit any towers, by the Kazhdan--Margulis Theorem in Section~\ref{s:covolumes} above.

The first examples of towers of tree lattices are due to Bass--Kulkarni~\cite{BK}.  Generalizing these constructions,
Rosenberg~\cite{Ros} proved that if $X$ is a tree such that $\Aut(X)$ is nondiscrete and admits a uniform lattice, then
$\Aut(X)$ admits a tower of uniform lattices.  Carbone--Rosenberg~\cite{CR} considered nonuniform lattice towers in $\Aut(X)$ for $X$ a tree, showing
that, with one exception, if $\Aut(X)$ admits a nonuniform lattice then it admits a tower of nonuniform lattices.

In higher dimensions, for $X$ a right-angled building (see Section~\ref{ss:rab}) Thomas~\cite{Th2} constructed a tower of uniform and of nonuniform lattices.  Other higher-dimensional cases are open.  In particular, it is not known whether the automorphism groups of any Fuchsian buildings which are not right-angled (see Section~\ref{ss:hyp_buildings}) admit towers.

A finer version of Question~\ref{q:towers} is the following:

\begin{ques} Does $G$ admit a tower of \emph{homogeneous} lattices, that is, lattices acting transitively on cells of maximum dimension in $X$?
\end{ques}

For $X=T_{p,q}$ the $(p,q)$--biregular tree, if $p$ or $q$ is composite there is a homogeneous tower in $G=\Aut(X)$ (Bass--Kulkarni~\cite{BK}).  When $X$ is the $3$--regular tree, a deep theorem of Goldschmidt~\cite{Go} implies that $G$ does not admit such a tower, since $G$ contains only finitely many conjugacy classes of edge-transitive lattices. The Goldschmidt--Sims conjecture (see~\cite{Gl}), which remains open, is that if $p$ and $q$ are both prime, then there are only finitely many conjugacy classes of homogeneous lattices in $\Aut(T_{p,q})$.  If $X$ is the product of two trees of prime valence,  Glasner~\cite{Gl} has shown that there are only finitely many conjugacy classes of (irreducible) homogeneous lattices in $G=\Aut(X)$.  For all other higher-dimensional $X$, the question is open.

\begin{ques} Does $G$ admit maximal lattices?
\end{ques}

In the algebraic setting, lattices of minimal covolume are known in many cases (see~\cite{Lu} and its references), and so these lattices are maximal.  Examples of maximal lattices in $G=\Aut(X)$ are some of the edge-transitive lattices for $X$ the $3$--regular tree, classified by Goldschmidt~\cite{Go}.

A coarse version of the question of towers is:

\begin{ques}[Lubotzky] Let $\G$ be a uniform lattice in $G=\Aut(X)$. Define
\[
   u_\G(n) = \# \bigl\{ \,\G' \bigm| \mbox{$ \G'$ is a lattice containing $\G$,
   and $[\G':\G] = n$}\bigr\}
\]
By similar arguments to~\cite{BK}, $u_\G(n)$ is finite.
What are the asymptotics of $u_\G(n)$? \end{ques}

The case $X$ a tree was treated by Lim~\cite{L1}.  If $X$ is the
building associated to a higher-rank algebraic group, then for any
$\G$, we have $u_\G(n)=0$ for $n>\!>0$, since $\mathcal{V}(G)$ has
positive lower bound.  In contrast, if $\Aut(X)$ admits a tower of
lattices (for example if $X$ is a right-angled building), there is a
$\G$ with $u_\G(n) > 0$ for arbitrarily large $n$.
Lim--Thomas~\cite{LT}, by counting coverings of complexes of groups,
found an upper bound on $u_\G(n)$ for very general $X$, and a lower
bound for certain right-angled buildings $X$.  It would be
interesting to sharpen these bounds for particular cases.

\subsection{Biautomaticity of lattices}
\label{ss:Biautomaticity}

The theory of automatic and biautomatic groups is
closely related to nonpositive
curvature.
All word hyperbolic groups are biautomatic \cite{ECHLPT}.
Yet it is not known whether an arbitrary group
acting properly, cocompactly, and isometrically on a CAT$(0)$
space is biautomatic, or even automatic.
Indeed, the following special case is open:

\begin{ques} Suppose a group $\G$ acts properly, cocompactly and isometrically on a CAT(0) piecewise
Euclidean $2$--complex.  Is $\G$ biautomatic? Is $\G$ automatic? \end{ques}

Biautomaticity is known in several cases for groups acting on
complexes built out of restricted shapes of cells.
Gersten--Short established biautomaticity for uniform
lattices in CAT$(0)$ $2$--complexes of type
$\tilde{A}_1 \times \tilde{A}_1$, $\tilde{A}_2$, $\tilde{B}_2$,
and $\tilde{G}_2$ in \cite{GerstenShort1,GerstenShort2}.
In particular, Gersten--Short's work includes CAT$(0)$
square complexes, $2$--dimensional systolic complexes,
and $2$--dimensional Euclidean buildings.

Several special cases of Gersten--Short's theorem have been extended.
For instance \Swiatkowski\ proved that any uniform lattice in a Euclidean
building is biautomatic \cite{Swiatkowski06_biautomatic}.
Niblo--Reeves \cite{NibloReeves98}
proved biautomaticity of all uniform lattices acting on
CAT$(0)$ cubical complexes.
In particular, this result includes all finitely generated
right-angled Coxeter groups and right-angled Artin groups.
Systolic groups, that is, uniform lattices acting on arbitrary systolic
simplicial complexes, are also biautomatic by work of
Januszkiewicz--\Swiatkowski~\cite{JS}.

Gersten--Short's work applies only to $2$--complexes with a single
shape of $2$--cell.  Levitt has generalized
Gersten--Short's theorem to prove
biautomaticity of any uniform lattice acting on a CAT$(0)$
triangle-square complex, that is, a $2$--complex
each of whose $2$--cells is either a
square or an equilateral triangle \cite{Levitt_TriangleSquare}.

Epstein proved that all nonuniform lattices in $\SO(n,1)$ are biautomatic
\cite[11.4.1]{ECHLPT}.
Rebbechi \cite{Rebbechi01} showed more generally
that a relatively hyperbolic
group is biautomatic if its peripheral subgroups are biautomatic.
Finitely generated virtually abelian groups are biautomatic
by \cite[\S 4.2]{ECHLPT}.
It follows from work of Hruska--Kleiner \cite{HK}
that any uniform lattice acting on a CAT$(0)$ space with isolated flats
is biautomatic.

By a theorem of Brink--Howlett, all finitely generated Coxeter groups
are automatic \cite{BrinkHowlett93}.
Biautomaticity has been considerably harder to establish, and remains
unknown for arbitrary Coxeter groups.
Biautomatic structures exist when the Coxeter group is affine,
that is, virtually abelian, and also when the Coxeter group has no affine
parabolic subgroup of rank at least three by a result
of Caprace--M\"uhlherr~\cite{CapraceMuhlherr05}.
Coxeter groups whose Davis--Moussong complex has isolated flats
are also biautomatic by \cite{HK}.
The Coxeter groups with isolated flats have been classified
by Caprace~\cite{Caprace_Isolated}.

Let $W$ be a Coxeter group, and let $X$ be a building of type $W$.
\Swiatkowski~\cite{Swiatkowski06_biautomatic}
has shown that any uniform lattice $\Gamma $ in $G=\Aut(X)$
is automatic.  If $W$ has a geodesic biautomatic structure,
he shows that $\Gamma$ is biautomatic as well.
Together with Caprace's work mentioned above
it follows that if $W$ is a Coxeter group with isolated flats,
then $\Gamma$ is biautomatic \cite{Caprace_Isolated}.
This consequence can be seen in two ways: using the fact that
$W$ is biautomatic, or alternately using the fact, established by
Caprace, that $\Gamma$ is relatively hyperbolic
with respect to uniform lattices in Euclidean buildings.

\pagebreak

\end{document}